\documentclass[usletter,12pt]{article}

\usepackage{latexsym,amssymb}
\usepackage[T1]{fontenc}
\usepackage[dvips]{graphicx}

\newcommand{\y}{\textquotedbl}
\newcommand{\s}{\textbullet}
\newcommand{\z}{\mathbb Z}

\newtheorem{lem}{Lemma}[section]

\newtheorem{example}{Example}[section]

\newtheorem{theorem}[lem]{Theorem}
\newtheorem{prop}[lem]{Proposition}

\newenvironment{proof}{\textbf{Proof.}}{\newline\hspace*{\fill}{$\Box$}}
\begin{document}
\title{Fibred and Virtually Fibred hyperbolic 3-manifolds in the censuses}
\author{J.\,O.\,Button\\
Selwyn College\\
University of Cambridge\\
Cambridge CB3 9DQ\\
U.K.\\
\texttt{jb128@dpmms.cam.ac.uk}}
\date{}
\maketitle
\begin{abstract}
Following on from work of Dunfield, we determine the fibred status
of all the unknown hyperbolic 3-manifolds in the cusped census.
We then find all the fibred hyperbolic 3-manifolds in the closed census 
and use this to find over 100 examples each of closed and cusped virtually
fibred non-fibred census 3-manifolds, including the Weeks manifold.
We also show that the co-rank of the fundamental group of every
3-manifold in the cusped and in the closed census is 0 or 1. 
\end{abstract}
\section{Introduction}

A famous open question of Thurston asks if every finite volume hyperbolic
3-manifold is virtually fibred, that is it
has a finite cover that is fibred over the circle. A finite
volume hyperbolic 3-manifold (which we assume throughout to be orientable)
is either closed or is the interior of a compact
3-manifold with boundary a finite union of tori, which we call the cusps. Let
us treat this as two separate questions, one about closed and one about cusped
3-manifolds. A reason put forward (for instance in \cite{kap}, \cite{lack})
as to why this question may not be true is that there are very few examples
known of non-fibred hyperbolic 3-manifolds that are virtually fibred. 
However we have data available in the form of the Callahan-Hildebrand-Weeks
census of nearly 5,000 cusped hyperbolic 3-manifolds and the
Hodgson-Weeks census of nearly 11,000 closed hyperbolic 3-manifolds which
should make a good testing ground. Computer programs run by Dunfield 
\cite{dunl}
show that over 87\% of the 3-manifolds in the cusped list are fibred, 
suggesting that non-fibred virtually fibred cusped hyperbolic 3-manifolds are 
not so easy to come by because fibred examples are so common.

This of course would not apply to closed 3-manifolds $M$ as if $M$ has
finite homology then it is not fibred, and this is the case for nearly all
3-manifolds in the closed census (although recently \cite{dunth} showed
with mammoth computation that they all have a finite cover with positive
first Betti number). 
In this paper we will find over 100 examples in the closed census of non-fibred
virtually fibred 3-manifolds, including 10 from the 30 with smallest
volume. All these examples are arithmetic and the first is the Weeks 
manifold, which is the one of minimum volume in the
census and conjectured to be the minimum volume hyperbolic 3-manifold
overall. Also one of the non-fibred virtually fibred 
examples has positive first Betti number, which is the first known case of
such a closed 3-manifold.

In order to do this we determine the fibred 3-manifolds in the cusped and
closed censuses. Our starting point is the list of Dunfield \cite{dunl} 
which used two programs to work out the fibred and non-fibred 3-manifolds 
in the cusped census,
with 169 exceptions which were left as unknown. We find the fibred
status of all of these unknowns: in fact 5 are fibred and 164 are not. 
After this we examine the 128 3-manifolds with positive first Betti number 
in the closed census and prove that 87 are fibred
with 41 that are not, thus providing the complete list of closed fibred
3-manifolds in the census. We then utilise the data given in the program
Snap and recent work of Goodman, Heard and Hodgson
to find other hyperbolic 3-manifolds which are commensurable with these
fibred ones, so are virtually fibred. 

All our
techniques only require knowledge of the fundamental group of the 3-manifolds,
as we can utilise a result \cite{stll} of Stallings. In particular we can
apply the Bieri-Neumann-Strebel (BNS) invariant and the Alexander polynomial to
these fundamental groups. In Section 2 we give a brief
description of the BNS invariant and demonstrate how it can sometimes be used 
to determine the fibred status of a hyperbolic 3-manifold, using a result
of K.\,S.\,Brown. We summarise the Alexander polynomial in Section 3.

In Section 4 we examine the unknown cusped 3-manifolds, by first applying the 
BNS invariant and the Alexander polynomial and then working
directly with the fundamental group. Next in Section 5 we use this
information and the knowledge of commensurability classes of cusped hyperbolic
3-manifolds to find non-fibred virtually fibred cusped
hyperbolic 3-manifolds. In Section 6 we obtain closed census fibred hyperbolic
3-manifolds from cusped ones. We do not quite 
pick up all closed fibred 3-manifolds from the census in this way, so then
we use the Alexander polynomial to demonstrate that most of the rest of the 
3-manifolds in the closed census with positive first Betti number are not
fibred, with those that remain shown to be fibred directly, using finite
covers. In Section 7 we then obtain closed non-fibred virtually fibred
hyperbolic 3-manifolds which are all arithmetic.

The co-rank of a finitely generated group is the largest integer $n$ for which
the group has a homomorphism onto the free group of rank $n$. To finish we
quickly show in Section 8 that all closed and cusped census 3-manifolds
have co-rank 0 or 1.

In the Appendix we have five tables: the first has the Alexander polynomials
of the unknown cusped census 3-manifolds and the second gives cusped
non-fibred virtually fibred hyperbolic census 3-manifolds.
The third displays all the closed fibred census 3-manifolds. Table 4 
lists all remaining
closed census 3-manifolds with positive first Betti number, so these are
exactly the non-fibred 3-manifolds in the closed census with positive 
Betti number, and Table 5 contains the closed non-fibred virtually fibred 
census 3-manifolds that we found.

We are taking as our input data
the two censuses which come with SnapPea, the related data in Snap
and with \cite{ghh}, the
presentations of fundamental groups from SnapPea as given in \cite{dunthda}
and the list \cite{dunl} of fibred 3-manifolds in the cusped census. 
From then on, we only work with a fundamental group
presentation and operate either by hand or by using a program that can
determine, and provide presentations for,
all subgroups of a given small index of a finitely presented
group, such as Magma or Gap. 
We would like to thank Craig Hodgson for introducing us to the censuses
and the referee for providing helpful comments and useful references
on receipt of an earlier draft of this paper.
\section{The Bieri-Neumann-Strebel Invariant}
If $G$ is a finitely generated group with $G'$ the commutator subgroup then
let $\beta_1(G)$ be the first Betti number of $G$, that is the number of
free summands in the abelianisation $\overline{G}=G/G'$. Assuming that
$b=\beta_1(G)>0$, there exist homomorphisms of $G$ onto $\mathbb Z$ and the
BNS invariant gives us information on when 
their kernels are finitely
generated. This is done in \cite{bns} by identifying non-zero homomorphisms
of $G$ into $\mathbb R$, up to multiplication by a positive constant, with
the sphere $S^{b-1}$. The BNS invariant of $G$ is an open subset $\Sigma$ of 
$S^{b-1}$, with a homomorphism $\chi$ of $G$ onto $\mathbb Z$ having finitely
generated kernel if and only if $\chi$ is in both $\Sigma$ and $-\Sigma$. If
$G=\pi_1M$ for $M$ the fundamental group of a compact 3-manifold then it is
shown that $\Sigma=-\Sigma$. In general it can be difficult to find
$\Sigma$ but in a paper of K.\,S.\,Brown \cite{ksb}, an
algorithm is given to determine whether or not $\chi$ is in $\Sigma$ in the 
case where $G$ is a one relator group. If $G$ has at least three generators
then $\Sigma=\emptyset$ so the interesting case is when we have a 2-generator,
1-relator group. But compact orientable irreducible
3-manifolds with non-empty toroidal boundary
always have a presentation with one less relator than the number of
generators and in the cusped census of 3-manifolds 
many (over 4000 out of 4815)
have 2-generator 1-relator fundamental groups.

The connection with fibred 3-manifolds dates back to a theorem of Stallings
\cite{stll} which states that if $M$ is compact, orientable and irreducible
with $\pi_1M$ possessing a surjection to $\mathbb Z$ with finitely
generated kernel then $M$ is fibred over the circle with the kernel being
the fundamental group of the fibre. Conversely if
$M$ is compact, orientable and fibred 
then of course $\pi_1M$ has this property 
and $M$ will be irreducible
except for $S^2\times S^1$: in fact as \cite{hemp} Chapter 11 makes clear,
if irreducibility is removed from the hypothesis of Stallings' result then
the conclusion still holds provided that $M$ has no sphere boundary
components (which we could cap off) and no fake 3-cells (for which we
could invoke the Poincar\'e conjecture). In any case we are interested in
hyperbolic 3-manifolds and these are always irreducible.

Thus the Brown algorithm will determine whether or not most 3-manifold in the
cusped census fibre. This is what Dunfield did, using a computer program to
work through the 3-manifolds $M$ which came with such a presentation and with
$\beta_1(M)=1$. The efficiency of the algorithm can be judged by the fact
that the total running time was about a minute. We outline how it works:
assume that $G=<a,b|r(a,b)>$ with $r$ reduced and cyclically
reduced. First suppose $\beta_1(G)=1$ so 
that there is one homomorphism $\chi$ from
$G$ onto $\mathbb Z$ (up to sign), with $\chi(a)=m$ and $\chi(b)=n$ (where
$m$ and $n$
can instantly be found by abelianising). Assume first that $m,n\neq 0$, then
we work through the relation, drawing a path which starts at height 0 and
rises or falls according to the value under $\chi$
of each successive letter in $r$. When
we finish, we must again be at height 0 and we regard this as being back at 
the starting point, having gone round in a circle. Then $\chi$ has finitely
generated kernel if and only if the path reaches both its maximum and its
minimum only once.

However one generator, say $a$, could have zero exponent sum which happens if
and only if $\chi(b)=0$, and then the criterion is slightly different: after
all there cannot now be a unique maximum. However in practice this case turns
out to be easier to work with, so we will make a definition: let us 
say throughout that a presentation
of a group $\Gamma=\langle g_1,\ldots , g_m|r_1,\ldots ,r_k\rangle$ with 
$\beta_1(\Gamma)=b\leq m$ is in standard form 
with respect to $g_1,\ldots ,g_b$ if
each of these has zero exponent sum in each relation $r_i$. Then these 
elements generate the infinite part of $\overline{\Gamma}$ with all other
generators being of finite order in $\overline{\Gamma}$. Now if 
$G=\langle a,b|r\rangle$ is in standard form, we have that
$\mbox{ker }\chi$ is finitely generated if
and only if the maximum and minimum occur twice, which will be either end of
a single flat path.

Given a compact orientable irreducible 3-manifold with $n$ cusps, we have by
Mayer-Vietoris that $\beta_1(M)\geq n$ so that this process can only work on
1-cusped 3-manifolds. But now suppose that our 2-generator 1-relator group
$G$ has $\beta_1(G)=2$. Then there are an infinite
number of homomorphisms from $G$ onto $\mathbb Z$ and here Brown's
algorithm works in the following way. 
We draw the (reduced and cyclically reduced) relation on a 2 dimensional
grid, and as it has zero exponential sum in both $a$ and $b$ we finish at the
origin. We then consider the convex hull $C$ in $\mathbb R^2$ of this path
and regard a homomorphism from $G$ onto $\mathbb Z$ as a directional vector,
with slope $n/m$ for $\chi(a)=m,\chi(b)=n$. Then the homomorphisms with
finitely generated kernel are those with slope lying between (but not
including) the slope of the outward pointing normals of
two successive edges of $C$, provided that the joining vertex, which will be a
vertex of the path, has only been passed through once when the path has been
traced out, along with the vertical homomorphism if and only if $C$ has a
unique horizontal side of length 1 on top, passed through only once, and 
similarly for the horizontal homomorphism. In fact a homomorphism is really
represented by two vectors with the same slope, pointing in opposite
directions, and both of these must satisfy the above conditions but again for
a 3-manifold group the conditions on each of the two vectors will be true
or false together because $C$ has rotational symmetry of order 2.
\begin{example}
\end{example}
\begin{figure}[h]
\begin{center}
\includegraphics[angle=-90,width=11cm]{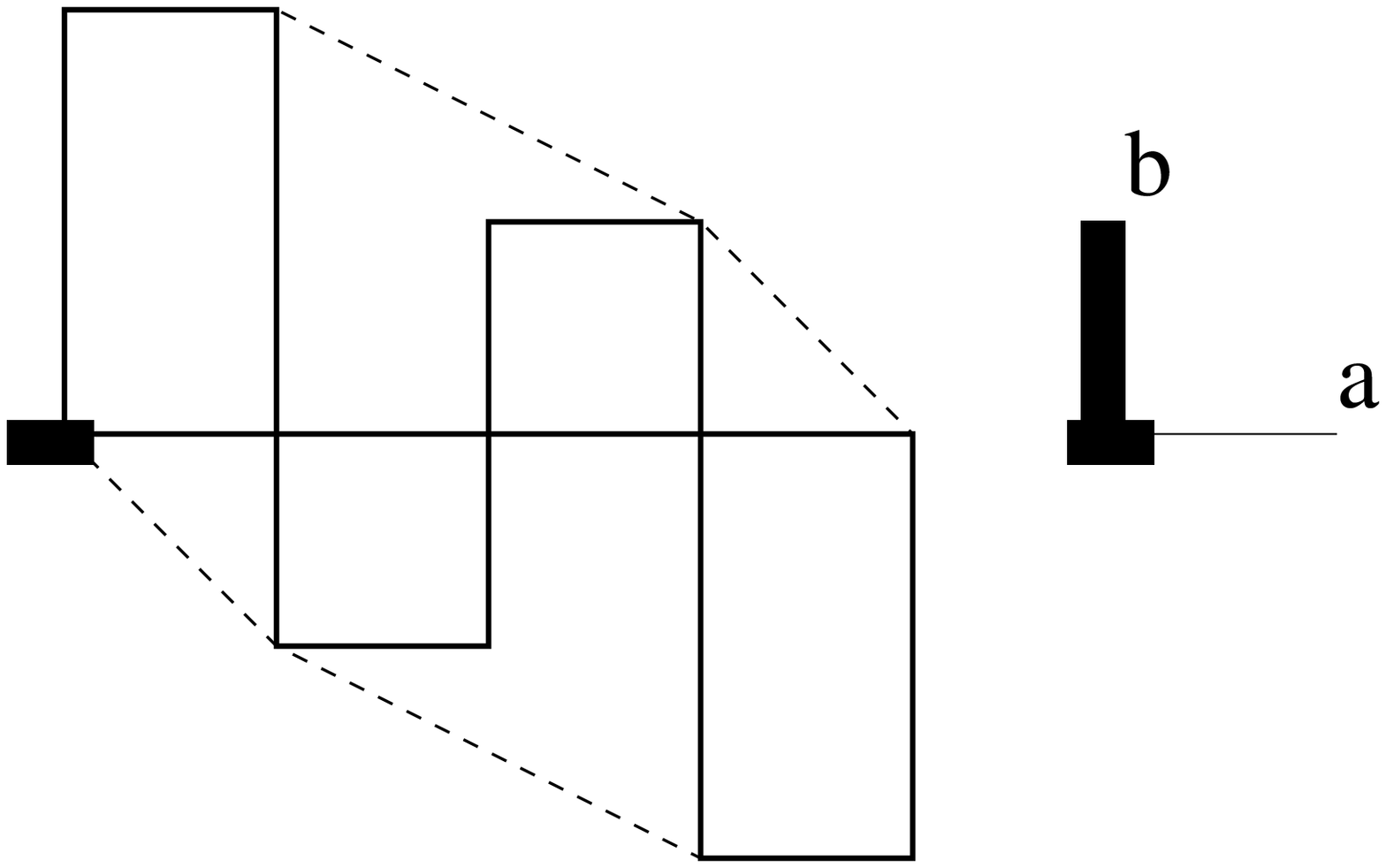}
\caption{}
\end{center}
\end{figure}
Let us demonstrate this process. We look for 1-cusped 3-manifolds in the 
census with $\beta_1(M)>1$ so that we have
a variety of homomorphisms to work with. We find only s789, v1539 and v3209,
all with homology $\mathbb Z+\mathbb Z$. SnapPea gives a 3 generator 
presentation for the fundamental group of two of them but we obtain
\[\pi_1(\mbox{v}1539)=\langle a,b|a^4B^2Ab^3AB^2Ab^3AB^2\rangle\]
with $(m,l)=(Ab,B^3a^5B^2)$ a basis for the fundamental group of the cusp.
This example will be important in Section 6. Drawing out the relation to
form the convex hull $C$ as in Figure 1 and using Brown's algorithm
reveal that all but the three homomorphisms (ignoring signs)
$\chi(a)=1,\chi(b)=0$; $\chi(a)=1,\chi(b)=1$ and $\chi(a)=1,\chi(b)=2$ have
finitely generated kernel.

Thus we see that determining the fibred status of cusped hyperbolic 3-manifolds
with a 2-generator 1-relator fundamental group presentation presents no
problem, but for a closed orientable irreducible 3-manifold $M$ we have that
every presentation of $\pi_1M$ has at least as many generators as relators.
Thus it would appear here that Brown's algorithm is now no use, however 
we make an obvious yet useful point: suppose we have
a 2 generator group $G=\langle a,b|r_1,\ldots ,r_m\rangle$ then any 
2 generator group $\Gamma$
of the form $\langle a,b|r\rangle$ where $r$ is one of the
$r_i$ (or even just a consequence of $r_1,\ldots ,r_m$) surjects onto $G$.
If we have a finitely generated kernel $K$ of a homomorphism $\chi$ from 
$\Gamma$ onto $\mathbb Z$, which can be determined by Brown's algorithm, 
then the image of $K$ in $G$ is still finitely generated, 
so the only issue is whether $\chi$ factors through $G$ and this is easily
solved by looking at the abelianisations of $\Gamma$ and $G$. In particular
if we have a surjection from any $\Gamma=\pi_1M$ to any $G=\pi_1N$ where
$M$ and $N$ are both compact orientable irreducible 3-manifolds 
with $\beta_1(N)=\beta_1(M)$ then $M$ fibred implies that $N$ is too.

An obvious method to obtain fundamental group surjections from 3-manifolds
to other 3-manifolds is through the use of Dehn surgery,
where we attach a solid torus to a component of the boundary of a cusped
3-manifold $M$. 
If the cusp has generators $m$ and $l$ in $\pi_1M$ then $(p,q)$ Dehn
filling for coprime integers $p,q$ with $q\geq 0$ means that we attach the
curve $m^pl^q$ to the compressible curve in the solid torus, thus adding
this relation to $\pi_1M$ and reducing the number of cusps by one. If
we start with a 1-cusped hyperbolic 3-manifold $M$ with $\beta_1(M)=1$
then there will be a 
unique Dehn surgery forming a closed 3-manifold $N$ with $\beta_1(N)=1$
(we might call this curve the longitude, in analogy with a knot in $S^3$
where this is the only simple closed curve on the boundary homologous to 0) and
thus if $M$ is fibred and $N$ is irreducible then 
$N$ is fibred too as the relevant homomorphism 
$\chi:\pi_1M\rightarrow \mathbb Z$ factors through $N$. In fact here we do not
need to know that $N$ is irreducible, as seen by picturing this geometrically, 
because we are just performing Dehn filling along the boundary slope of the
fibre of $M$. This observation will be used
in Section 6, but to conclude this section let us apply this to
our example $M=\mbox {v1539}$. Performing $(p,q)$ Dehn surgery with the 
above basis
for the cusp means that the only homomorphism $\chi$ that factors through
$\pi_1N$ is $\chi(a)=\chi(b)=1$ (unless $(p,q)=(5,1)$ in which
case they all do) which is one of the three exceptional homomorphisms
so this does not tell us that $N$ is fibred.
However we can use the Dehn filling relation instead to give us:
\begin{figure}[h]
\begin{center}
\includegraphics[angle=-90,width=9cm]{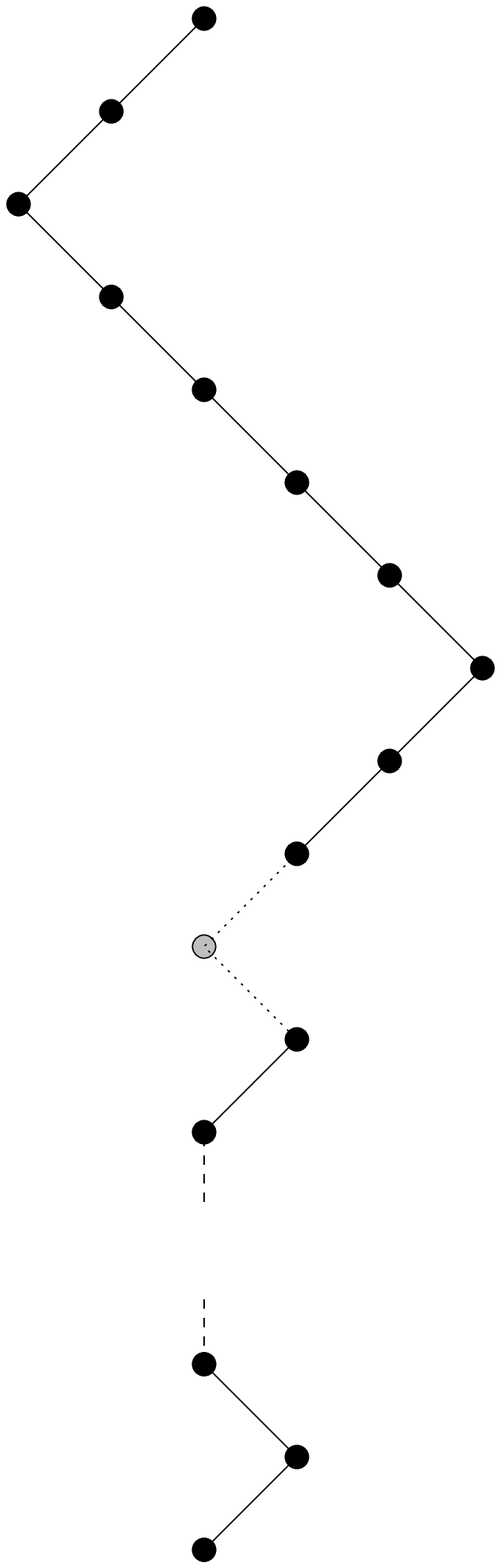}
\caption{}
\end{center}
\end{figure}
\begin{theorem} There exist infinitely many closed hyperbolic fibred 
2-generator 3-manifolds with bounded volume.
\end{theorem}
\begin{proof}
We take v1539($p$,1) and consider $\Gamma=
\langle a,b| m^pl\rangle$ which surjects 
onto its fundamental group, with $\beta_1(\Gamma)$ 
also equal to 1 if $p\neq 5$.
Taking the homomorphism $\chi(a)=\chi(b)=1$, we draw out the relation
as in Figure 2,
where we have cancellation along the dotted lines if $p>0$ but we still have
a unique maximum and minimum, hence a finitely generated kernel. 
We then apply Thurston's Dehn surgery theorem to obtain
hyperbolicity, hence irreducibility which gives us  
the fibred property, along with the fact that these
closed 3-manifolds have volume accumulating to that of v1539.
\end{proof}
\section{The Alexander polynomial}
Historically the Alexander polynomial was first introduced for knots in 
$S^3$ but it can be defined for any finitely presented group.
Although it is not able to give us so much information as the BNS invariant,
it has the advantage that it is straightforward to work
out from any finite presentation of a group using Fox's free differential
calculus. Therefore we give a brief description adopting the approach of
Fox in \cite{fx}.

Let the finitely presented group $G$ be $\langle x_1,\ldots ,x_n|r_1,\ldots ,
r_m\rangle$ in terms of generators and relators, and let its free 
abelianisation be $ab(G)$, which will be isomorphic to 
$\mathbb Z^b$ where $b=\beta_1(G)$.
If $F_n$ is the free group of rank $n$ with free
basis $x_1,\ldots ,x_n$ then a derivation of the integral group ring
$\mathbb Z[F_n]$ is a map from $\mathbb Z[F_n]$ to itself satisfying
\begin{eqnarray*}D(v_1+v_2)&=&Dv_1+Dv_2,\\ 
 D(v_1v_2)&=&(Dv_1)\tau(v_2)+v_1Dv_2
\end{eqnarray*}
where $\tau$ is the trivialiser: namely the ring homomorphism from 
$\mathbb Z[F_n]$ to $\mathbb Z$ with $\tau(x)=1$ for all $x\in F_n$. It is a
fact that for each free generator $x_j$ there exists a unique derivation
$D_j$, also written $\partial/\partial x_j$, such that $\partial x_i/\partial
x_j=\delta_{ij}$. 
To calculate the ``partial derivative'' $\partial w/\partial x_j$ for any
$w\in F_n$ we can use the formal rules
\[\frac{\partial x_i}{\partial x_j}=\delta_{ij},\qquad
\frac{\partial x_i^{-1}}{\partial x_j}=-\delta_{ij}x_i^{-1},\qquad
\frac{\partial (w_1w_2)}{\partial x_j}=\frac{\partial w_1}{\partial x_j}+w_1
\frac{\partial w_2}{\partial x_j}\]
where generally $w_2$ will be the last letter in the word $w=w_1w_2$. Let
$\gamma$ be the natural map from $\mathbb Z[F_n]$ to $\mathbb Z[G]$ and let
$\alpha$ be the same from $\mathbb Z[G]$ to $\mathbb Z[ab(G)]$. Then the
Alexander matrix $A$ of the presentation is the $m\times n$ matrix with
entries
\[a_{ij}=\alpha\gamma\left(\frac{\partial r_i}{\partial x_j}\right).\]
We define the $k$th elementary ideal $E_k(A)$ to be the ideal of
$\mathbb Z[ab(G)]$ generated by the $(n-k)\times (n-k)$ minors of $A$
if $0<n-k\leq m$, thus
under this notation $k$ is the number of columns that are deleted in forming
the minors. Finally we define the Alexander polynomial 
$\Delta_G$ to be the generator (up to units) of the
smallest principal ideal containing $E_1(A)$. To calculate it 
we can choose a basis $(t_1,\ldots ,t_b)$
for $ab(G)$, apply the free
differential calculus as above and then form our matrix by evaluating. 
From here we can determine the minors and their highest common factor.
Of course this would be of little use if it depended on the presentation
of $G$, but that it is invariant can be seen directly, as shown in
\cite{fx} VII 4.5, by observing that applying a Tietze transformation to a
presentation does not change the elementary ideals. Alternatively we have
a topological definition of the Alexander polynomial, as described 
in \cite{mcm} Section 2 or \cite{dunp} Section 3: if $X$ is a finite CW-complex
with $\pi_1X=G$ and $f:\tilde{X}\rightarrow X$ is the regular cover 
corresponding to the homomorphism $\alpha$ from $G$ to 
$ab(G)$ then, taking $p\in X$, the Alexander
module of $X$ over the group ring $\mathbb Z[\mbox{ab}(G)]$ is
$H_1(\tilde{X},f^{-1}(p);\mathbb Z)$. The connection between the two approaches
is that by taking a free resolution of this module, we obtain the Alexander
matrix as above (or rather under our notation it is the transpose of $A$).
The Alexander polynomial $\Delta_G$ is  only defined up to
units, thus we can think of
$\Delta_G$ as a Laurent polynomial in $\mathbb 
Z[t_1^{\pm 1},\ldots ,t_b^{\pm 1}]$ up to multiplication by 
$\pm t_1^{k_1}\ldots t_b^{k_b}$.
Of course the actual coefficients depend on this basis: sometimes there will
be a natural choice, such as for a $b$-component link in $S^3$ where we
would take meridians about each link. However we might not in general have
this luxury, although we can always
make a change of basis if necessary by putting
$t_i=s_1^{k_{i1}}\ldots s_b^{k_{ib}}$ with the vectors 
$(k_{i1},\ldots ,k_{ib})$ making up an element of $GL(b,\mathbb Z)$.

The utility of the Alexander polynomial for us here is the well known result,
derived later,
that if we have a compact 3-manifold $M$ with $\beta_1(M)=1$ then its
Alexander polynomial $\Delta_M(t)$, in this case a Laurent polynomial
defined up to units and with $\Delta_M(1/t)$ equal to $\Delta_M(t)$ times a
unit, is monic if $M$ is fibred. We also have by Dunfield a suitable
generalisation of this for the case $\beta_1(M)\geq 2$ which we will use later:
Theorem 5.1 of \cite{dunp} states that if the Alexander
polynomial $\Delta_M$ has no terms with coefficients that are $\pm 1$ then
$M$ is not fibred: more precisely let $N$ be the Newton polytope of 
$\Delta_M$, that is the convex hull in $\mathbb R^b$ of the points 
$(k_1,\ldots ,k_b)$ where $x_1^{k_1}\ldots x_b^{k_b}$ is a (non-trivial)
term of $\Delta_M$. If none of the vertices of $N$ have coefficient $\pm 1$
in $\Delta_M$ then the Bieri-Neumann-Strebel invariant $\Sigma$ of $\pi_1M$
is empty and so there are no homomorphisms onto $\mathbb Z$ with finitely
generated kernel.

\section{The unknown cusped 3-manifolds}

When Dunfield ran his programs on the 4815 3-manifolds in the cusped census to
see which were fibred, he first
set up the computer to apply Brown's algorithm to
any 3-manifold $M$ with a 2 generator 1 relator presentation and with
$\beta_1(M)=1$. As we have seen in Section 2, this is guaranteed to 
terminate and give a definite 
yes/no answer. The program took about a minute in total to complete the 4105
examples given to it, 3653 of which were fibred and 452 of which were not.

The other algorithm that was applied was Lackenby's idea of
taut ideal triangulations. We will not be using this because our emphasis is
on methods which only require knowledge of the fundamental group; we note 
only that this process will not tell us that the 3-manifold is non-fibred
but it has no restriction as above on the number of generators or relators.
When this was applied to the cusped census it produced 541 further fibred
3-manifolds, as well as confirming a lot of the 3-manifolds already known
to be fibred by Brown's algorithm. There were some of these that it did not
work for, and the running time was a lot longer.

Thus this leaves 169 cusped 3-manifolds whose status is unknown. In
this section we will determine whether or not these are fibred.
As any unknown 3-manifold has already passed
through the two algorithms above, we proceed by a variety of listed
methods involving fewer and fewer
3-manifolds. We work on the assumption that they are most likely to be
non-fibred, because a fibred 3-manifold has had two chances already to be
detected, and then only at the very end do we admit the possibility that
what remains might be fibred.

{\bf 1. Use other data}\\
In \cite{cdw}, all knots in $S^3$ appearing in the m or s part of the census
are determined and listed, helpfully with the genus of their fibre or an
{\texttt x} if they are non-fibred. We might as well annotate
Dunfield's list to provide a fuller description of such 1-cusped 3-manifolds.
We find ourselves marking an unknown 3-manifold
on 3 occasions: m372 is the non-fibred knot $9_{46}$ in the Alexander-Briggs/
Rolfsen-Bailey tables (for alternative names we have 3,3,2 1- in Conway
notation or $9n5$ in the Dowker-Thistlethwaite ordering used in Knotscape,
where $n$ denotes a non-alternating knot), s879 is a non-fibred
knot with 11 crossings (5,3,2 1- or $11n139$),
and s704 is the fibred knot $10_{140}$ (equivalently 4,3,2 1-
or $10n29$) with genus 2.
(This is somewhat lucky - very few of the remaining 3-manifolds are fibred).

{\bf 2. Any other 2 generator groups?}\\
In the course of our study, we found one 3-manifold $M$ with a 2-generator
1-relator presentation and with $\beta_1(M)=1$ which was listed as unknown.
This is v3036 with presentation
\[a^3b^3AbAb^3a^3b^3AbAb^4AbAb^3\]
which we see is in standard form with respect to $a$. On applying Brown's
algorithm,
we reach the top after the middle $a^3$ term whence we have $b^3$, 
so this is not fibred.

We also find two 2-generator 1-relator 3-manifolds $M$ with $\beta_1(M)=2$
and with status unknown,
for which we can use the extended version of Brown's algorithm. We can
quickly check these are all the unknowns of this form because
the cusped census collects 3-manifolds with the same number of cusps
together. But $\beta_1(M)$ is at least the number of cusps and we know that
there are only three cases 
where $M$ has one cusp but $\beta_1(M)=2$, with these listed as
fibred. Therefore we work down the table of 2-cusped 3-manifolds, all of
which happen to have $\beta_1(M)=2$, and look them up in Dunfield's list. We
know that either they will be proved fibred using taut foliations or they
will be unknown. In fact we find that it is the former in all but four cases:
v2943, v3379, v3384, v3396. The last two have homology $\mathbb Z_5+\mathbb Z
+\mathbb Z$ and $\mathbb Z_3+\mathbb Z+\mathbb Z$ respectively so are not
2 generator, but we find
\begin{eqnarray*}\pi_1(\mbox{v}2943)&=&\langle a,b|
abAB^2AbaBAba^3bABab^2aBAbaBA^3B\rangle,\\
\pi_1(\mbox{v}3379)&=&
\langle a,b|abABa^3BAbaBAbaB^2abABabA^3baBAbaBAb^2AB\rangle,
\end{eqnarray*}
neither of which are fibred, seen by drawing out the relation and noting that
all vertices of the convex hull are passed through more than once.

Moreover there are only three 3-cusped 3-manifolds $M$, all of which are
fibred and have $\beta_1(M)=3$, and none at all with more than three cusps.
This now leaves only 1-cusped 3-manifolds, apart from v3384 and v3396. 

{\bf 3. The Alexander Polynomial}\\
We now turn to the the original suggestion of Dunfield of
calculating Alexander polynomials. Once some practice is gained, the process
becomes much faster so we might as well apply it to all the remaining
unknowns. Let us first assume that $M$ is a 
1-cusped 3-manifold with $\beta_1(M)=1$. As mentioned in Section 3, on taking
$t$ as a generator (by symmetry it does not matter which one) for $ab(\pi_1M)$
we have that the Alexander polynomial of $M$ is an element of the ring 
$\mathbb Z[t,t^{-1}]$, up
to units which are $t^{\pm k}$ for $k\in\mathbb Z$. 

In the process of calculating the polynomial, we found it quickest to make
substitutions so that we always have a presentation for
$\pi_1M$ which is in standard form with respect to
one of the generators, say $x$. Then it is seen that 
$\partial r_i/\partial x=0$ on evaluation for each of the relations $r_i$:
first note that $\alpha(g_j)=1$ for all the other generators $g_j$ of our 
presentation. Thus whenever we have an $x$ appearing in
$r_i$ it contributes a term which is (on evaluation) $t^k$, where $k$ is the
exponent sum of $x$ in the subword of $r_i$ strictly to the left of this
appearance of $x$, whereas an $X$ contributes $-t^k$ for $k$ the exponent
sum of $x$ in the subword to the left of and including $X$. The result then
follows by pairing off each $x$ and the $X$ with which it cancels when all
other $g_j$ are set to the identity. A special case of a presentation in
standard form is when each relator
has only one appearance of $x$, which we refer to as simple form with respect
to $x$, so we get
\begin{equation}
r_i=xu_iXv_i\qquad\mbox{and}\qquad
\frac{\partial r_i}{\partial g_j}=k_{ij}t+l_{ij}
\end{equation}
where $u_i,v_i$ contain no appearance of $x$ and $X$, with $k_{ij}$ the
exponent sum of $g_j$ in $u_i$ and $l_{ij}$ that of $g_j$ in $v_i$. In
particular if $M$ is fibred over the circle with fibre the surface $S$, so
that $\pi_1S$ is free of rank $n$, then we can take a presentation for
$\pi_1M$ of the form $\langle g_1,\ldots ,g_n,x|r_1,\ldots ,r_n\rangle$,
where $r_i=xg_iXv_i$. Thus 
$\partial r_i/\partial g_j=\delta_{ij} t +l_{ij}$ so that the Alexander
polynomial is the characteristic polynomial of the $n\times n$
monodromy matrix $-l_{ij}$
induced by the glueing map, and hence is monic with degree $n$. Thus we look
for non-monic Alexander polynomials in our calculations and conclude that
these 3-manifolds are non-fibred.

In fact in the case of a 2-generator 1-relator group $G$ with $\beta_1(G)=1$
there is a straightforward connection between Brown's algorithm
and the Alexander polynomial $\Delta_G$: the way to see this is to assume that
$G=\langle a,b|r\rangle$ is in standard form with respect to $a$
and then once the relation is drawn out we note that the
process given of calculating $\Delta_G$ is merely that of counting the
appearance of $b$s (which contribute $+1$) and $B$s ($-1$) in the relation
at each level, and these values are the coefficients of $\Delta_G$.
In particular we obtain a very visual insight into how a 2-generator 
1-relator knot could have monic Alexander
polynomial but not be fibred; the relation must reach its peak more than
once but all but one of them must cancel out. Another example is that we
can easily recognise 1-punctured torus bundles amongst hyperbolic
3-manifolds with 2-generator 1-relator fundamental groups;
if $\pi_1M=\langle a,b|r\rangle$ with $r$ reduced and cyclically reduced
is the fundamental group of a hyperbolic
3-manifold $M$ then $M$ is a 1-punctured torus bundle if and only if
$\beta_1(M)=1$ and the relation lies on only three levels
with a unique maximum and
minimum when drawn out in standard form. This is because
hyperbolic 1-punctured torus bundles $M$ must have $\beta_1(M)=1$ and the 
other condition is exactly what is needed to conclude that $M$ fibres with
Alexander polynomial of degree 2, thus the fibre must be a 1-punctured
torus or a 3-punctured sphere, but the bundle is not hyperbolic in the
latter case. 
Now 1-punctured torus bundles might need three generators, as seen by
looking at their homology, but we cannot conclude in general that a 
hyperbolic 3-manifold $M$ is a 1-punctured torus bundle if it has a monic
quadratic Alexander polynomial. However, if we already know that $M$ is
fibred then we can.

Returning to the unknown cusped 3-manifolds $M$,
all our calculations are on 3 generator 2 relator groups so that we put
$\pi_1M=\langle g_1,g_2,x|r_1,r_2\rangle$ into standard form with
respect to $x$ and then we calculate the determinant of the $2\times 2$
matrix $\partial r_i/\partial g_j$. If furthermore 
our two relations are in simple form with respect to $x$,
that is as in (1) which happens often, 
then we can take a shortcut as the Alexander polynomial
will be (at most) quadratic. 
We calculate det$(k_{ij})$ which will be the coefficient
of $t^2$, and then det$(l_{ij})$ which is the constant. These must be equal
which acts as a useful check, given that we are doing these by hand (and
are here not interested in the middle term). 
More generally we ensure that our result is a Laurent polynomial that is 
symmetric under $t\mapsto t^{-1}$. The results are listed in Table 1 with
only six of these unknown 3-manifolds, written in bold, having a monic
Alexander polynomial. We can draw definite conclusions for two of them:
recall from Part 1 that s704 is a fibred knot, whereas
v2530 with Alexander polynomial $t+1$ of degree 1 cannot be fibred because the
fibre subgroup would have to be cyclic.

We can see from the table that some properties of the Alexander polynomial
of a knot are no longer true in this wider setting: for instance  
we no longer have $|\Delta_M(1)|=1$. In fact we can see from
our method of calculation of $\Delta_M$ on a presentation in standard form
that for $t=1$ we are just forming the equations of the
exponent sums of those generators (all but one) which have finite order in
homology, so $\Delta_M$ is never zero because
$|\Delta_M(1)|$ is always the order of the finite part of the homology.
(As this was not known to us when first compiling the table, it provided
another useful check). We can even have a common
factor of all the coefficients, as in $\Delta_{s773}(t)=2(t^3+t^2+t+1)$.
Moreover this example shows that Alexander polynomials are not necessarily of
even degree as they are for knots; other examples would be if $M$ is fibred
over a surface with an even number of boundary components (whereas knots can
only be fibred over a surface with one boundary component). 

We also need to consider the unknown 2-cusped 3-manifolds
v3384 and v3396. Taking the given presentation for $G=\pi_1(\mbox{v}3384)$ 
and putting it into standard form with respect to $(b,c)$ 
via the substitution $a=yB^2$ gives us the two relations
\[y^3B^2Cb^2y^2Bcb,\qquad yB^2Cb^2Yc\]
so the Alexander matrix is on evaluation (ordering the generators as $(b,c,y)$
and using the images of $b,c$ in $ab(G)$ as a basis, for which we also write
$b,c$):
\[\left(\begin{array}{rrr}
b^{-2}c^{-1}(1-c)&b^{-2}c^{-1}(b-1)&c^{-1}(2+3c)\\
b^{-2}c^{-1}(1+b)(1-c)&b^{-2}c^{-1}(b-1)(b+1)&c^{-1}(c-1)\end{array}\right)\]
giving the three minors (up to units):
\begin{eqnarray*}
m_1&=&-(b-1)(3bc+2b+2c+3)\\
m_2&=&(c-1)(3bc+2b+2c+3)\\
m_3&=&0\end{eqnarray*}
thus the Alexander polynomial is $3bc+2b+2c+3$. Similarly the given
presentation for $\pi_1(\mbox{v}3396)$ is already in standard form with respect
to $(b,c)$ so adopting the same notation we find its Alexander polynomial
is $2(b-1)(c-1)$. As mentioned at the end of Section 3, this gives us that
v3384 and v3396 are not fibred.

{\bf 4. Fibred after all?}\\
We now have to face up to the four remaining unknowns s594, v2869, v3093, 
v3541, and should take seriously the possibility that they are fibred.
If so then we must have a presentation
\begin{eqnarray}
\pi_1M=\langle t,a_1,\ldots ,a_r|ta_it^{-1}=w_i\rangle
\end{eqnarray}
where each $w_i$ is a word in $a_1,\ldots ,a_r$ equal to $\phi_*(a_i)$, for
$\phi_*$ the induced automorphism of $\pi_1M$ obtained from the glueing
homeomorphism $\phi$. These words, as well as $a_1,\ldots ,a_r$, 
generate the fibre
subgroup $F$ which will be free of rank $r$ equal to the degree of the
Alexander polynomial. Such a  presentation
will need more than the three generators that we have been given for our
3-manifolds, and
it might not be easy to move between the two different presentations.
However some points are clear:
as $\beta_1(M)=1$, the elements of $F$ are precisely those in $\pi_1M$ with
finite order in homology, and in looking for a candidate for $t$, any element
generating the infinite part of the homology can be used because we can
replace $t$ with $kt$ for any $k\in F$,
and $w_j$ with $kw_jk^{-1}$ in the presentation above.

In order to get round the number of generators, we use finite covers.
If $\pi_1M$ is fibred then we will have the cyclic covers $\pi_1M_n$ of
degree $n$, generated by the $r+1$ elements $t^n,a_1,\ldots ,a_r$ and with
$r$ relations, which correspond to the glueing homeomorphisms $\phi^n$.
When we ask Magma for a presentation of an index $n$ subgroup
of our 3 generator 2 relator group, it employs the Reidermeister-Schreier
process which will obtain a presentation of $2n+1$ generators and $2n$
relators, but some of these might be redundant so the output could be less.
Therefore we start with our unknown $\pi_1M$, using a presentation in standard
form with respect to a generator $x$. We ask Magma for (the generators of)
subgroups of index $n$ (it gives a subgroup in each conjugacy class)
and pick the cyclic cover $H_n$, that is the one with
the exponent sum of $x\equiv 0\mbox{ mod }n$ (which is easy to spot by
checking this condition holds for all of the given generators). 
We then demand a
presentation of $H_n$, hoping not only that it is $d+1$ generator and $d$ 
relator for $d$ the degree of the Alexander polynomial, but also that the
presentation $\langle h,x_1,\ldots ,x_d|r_1,\ldots ,r_d\rangle$ is in simple
form with respect to the generator $h=x^n$ of $H_n$. Then we look at the $d$
subwords from $h$ to $h^{-1}$ in each relation and
if this is a basis for the
free group on $x_1,\ldots ,x_d$ we conclude that conjugation by $h$ sends
$\langle x_1,\ldots ,x_d\rangle$ into itself. If now the subwords appearing
from $h^{-1}$ to $h$ are also a basis then $\langle x_1,\ldots ,x_d\rangle$
is normal in $H_n$, with $H_n$ having a presentation exactly as in (2) 
so by Stallings' condition
we have a finite cover of $M$ which is fibred, with fibre subgroup
$\langle x_1,\ldots ,x_d\rangle$. In fact we can halve the work as we need
only check that one of the two sets of subwords is a basis. This follows
from Proposition 3.1 in K.\,S.\,Brown's paper \cite{ksb}: suppose that $G$
is a finitely generated group and $\chi:G\rightarrow\mathbb Z$ is a
surjective homomorphism. To say that a HNN decomposition of $G$
has $\chi$ as associated homomorphism means that we can write $G$ as
$\langle B,t|B_1=tB_2t^{-1}\rangle$, for $B$ a subgroup of $G$ and $B_1,B_2$
subgroups of $B$, with $\chi(B)=0,\chi(t)=1$. Then we use the result that
$\chi\in\Sigma$ if and only if every HNN decomposition of $G$ with $\chi$ as
associated homomorphism is ascending, namely $B_2=B$. If this is so then we
can further ask whether $\chi\in -\Sigma$, but $-\chi$ is associated with
the decomposition of $G$ where $B_1$ and $B_2$ are swapped, thus a second
yes answer implies that $B_1=B_2=B$. However if $G$ is a 3-manifold group
then $\Sigma=-\Sigma$, meaning that one condition is enough.

To move from the fibred cover back to the original 3-manifold we use 
\cite{but} Corollary 2.6 which says 
that if the fibred 3-manifold $N$ is a finite cover of the compact orientable
3-manifold $M$, so that
$\beta_1(N)\geq\beta_1(M)$, then $M$ is fibred if the natural map given by
inclusion between the
infinite part of the abelianisations $\overline{\pi_1N}$ to 
$\overline{\pi_1M}$ has kernel coming from the fibre subgroup of $N$.
But if $\pi_1N$ is equal to $H_n$ as above and $x_1,\ldots ,x_d$ are elements
of finite order in the homology of $\pi_1M$ (which just means that when 
expressed as elements of $\pi_1M$ they have zero exponent sum in $t$) then,
as $h$ has infinite order in $\overline{\pi_1M}$, we have that the kernel
will be generated by $x_1,\ldots ,x_d$ (considered as elements of 
$\overline{\pi_1N}$) so it will be contained in the fibre subgroup of $N$. We
shall see directly that this condition always holds so we can conclude that
$M$ is fibred as well.

Starting with s594, the Alexander polynomial has degree 3 so, using the 
presentation $\langle a,c,x\rangle$ in standard form with respect to $a$
as obtained from Table 1, we
see that the index 2 subgroup $H$ corresponding to the cyclic cover has
abelianisation $\mathbb Z_2+\mathbb Z_4+\mathbb Z+\mathbb Z$, so is at least
four generator. On rewriting we are told it is generated by $p=x,q=c,r=axa^{-1}
,t=a^2$ with relations
\[
RQR{\bf t}pqp{\bf T},\quad PQP{\bf T}qP{\bf t}P,\quad QP{\bf T}R{\bf t}RqP
\]
and taking the subwords between $T$ and $t$ we easily see that these generate
the free group on $p,q,r$ so the cover is fibred, as is s594.
We can detect the fibre by noting that it must have fundamental group free
of rank 3, so is a 4-punctured sphere or a 2-punctured torus. In fact it
must be the latter because the glueing homeomorphism must permute the
boundary components and any one that is fixed must be sent to a conjugate of
itself in the fundamental group of the fibre under the induced automorphism
(it is not sent to
its inverse as the map is orientation preserving), thus adding 1 to
the Betti number of the 3-manifold. Thus if we have a 4-punctured sphere for
s594 then as it has Betti number 1, the induced permutation must be without
fixed points. But we can check that the cyclic cover of degree 4 has Betti 
number 3, whereas we would need the answer 5.

Moving onto v3093, we have
$\pi_1(\mbox{v}3093)=\langle b,x,y\rangle$ in standard form
with respect to $b$ and with degree four Alexander polynomial. Looking with
Magma at the finite index subgroups, the fundamental groups $H_n$ of the 
cyclic covers
of degree 2 and 3 are given with four generators, whereas of degree 4 and 5
we have 6 generators. On rewriting this cannot increase, so we try the
rewriting process for $H_4$ and $H_5$ 
which do then have the required 5 generators and
4 relations, with $t=b^n$ appearing as a generator. Unsurprisingly $t$ appears
too many times in the relations for $H_4$ but luckily we have $H_5$ with
abelianisation $\mathbb Z_2+\mathbb Z_2+\mathbb Z_2+\mathbb Z_2+\mathbb Z$
in simple form with respect to $t$:
setting $p=x,q=y,r=b^{-1}xb,s=b^{-1}yb$ and $t=b^5$ we have
relations
\[rsQPqrp^2q{\bf T}qSrsPRs{\bf t},\quad sP{\bf t}QPqrpqrpqr{\bf T}SrsPsP,\]
\[Qpq{\bf T}pSQpqSrsPRsPq{\bf t}RQPR,\quad
s{\bf T}QpqSrsPRsqSrsPRsqSrsPRs{\bf t}r^2,\]
and we get the computer to show that the subwords between $t$ and $T$
are a basis, by
setting up a homomorphism from the free group $F_4=\langle p,q,r,s\rangle$
to itself with these as images, and asking if it is a surjection. It is.
(We later confirmed this by hand, after obtaining practice with similar
calculations in Section 6.)

With the two remaining unknowns, v2869 and v3541, their Alexander polynomials
have degree 6 and 10. For v2869 we need a subgroup of at least index 3 to
have a hope of 7 generators, but the cyclic covers of degree 3,4,5 all fall
short. For v3541 we need index at least 5 for 11 generators, but index 5,6,7
all have 8 or less generators on rewriting. On trying to list all subgroups of
higher index we run into the problem that there are just too many. Instead
we rely on the fact that we have a good idea what the generators of these
particular cyclic covers should look like: if our original
fundamental group $G=\langle u,v,t\rangle$ 
is in standard form with respect to $t$
then $H_n$ has a generating set $t^iut^{-i},t^jvt^{-j},t^n$ for various values
of $i,j$, and on guessing such a generating set we can ask for the index of
$H_n$ in $G$ to check we are correct. Therefore, as we have 
$\pi_1\mbox{(v2869)}=\langle x,y,z\rangle$ in standard form with respect to 
$x$, we look at the subgroup $H$ 
generated by $x^iyx^{-i},x^jzx^{-j},x^n$ for
$i=0,\pm 1,\pm 2$, $j=0,\pm 1,-2$ and $n=6$. We do indeed find that $H$ has
index 6 in $G$ with abelianisation $\mathbb Z_{13}+\mathbb Z$ and on rewriting
we get the magic 7 generator 6 relator presentation, with generators
\[(a,b,c,d,e,f,t)=(y,z,xyX,Xzx,x^2yX^2,x^2zX^2,x^6)\]
which is in simple form with respect to $t$ and with the following free
basis to be found between $T$ and $t$:
\[
(F^2eBdBaceBabf,Fef,FBabFeBAbEceBabf,F^2eBdBaeBadBaCAbDbEf^2,\]
\[F^2eBdBacAbDbEFeBdBabFeBAbEcef,F^2eBdBacAbDbf).\]
Finally for $\pi_1\mbox{(v3541)}=\langle x,y,z\rangle$ in standard form with
respect to $z$, we try the subgroups $H_n$ generated by $z^ixz^{-i},z^jyz^{-j},
z^{-n}$ for $i=0,\pm 1,\pm 2$, $j=0, \pm 1, \pm 2$ and with $n$ running
from 8 to 15.
All have the correct index: for $n=8,9,10$ we get too few generators again on
rewriting but for the other $n$ we get exactly the required 11 generators
and 10 relations.
For $n=11$ the presentation is in standard but not in simple form with
respect to $t=z^n$, for the others it is indeed in simple form but with the
relations becoming progressively longer, so we take $n=12$. The subgroup
has abelianisation $\mathbb Z_7+\mathbb Z_{35}+\mathbb Z+\mathbb Z+\mathbb Z$
with the
other 10 generators
\[(a,b,c,d,e,f,g,h,i,j)=(x,y,zxZ,zyZ,Zxz,Zyz,z^2xZ^2,Z^2yz^2,Z^3yz^3,Z^4yz^4).
\]
Happily we find a basis between $t$ and $T$ of the form below:
\begin{eqnarray*}
(WJ,jwJ,jwfBAweJiCIjaI,jI,jEWaJ,iAJicI,jbDIjEWabFeWJ,&&\\
iH,hFEfBAweJidBAweJ,hCIhgFEfBAweJidFWJ)
\end{eqnarray*}
where $w=bDCIjaIhGHicH$. 

We have already mentioned in Part 1 the paper \cite{cdw} which lists the
knots in $S^3$ from the m and s part of the census. Recently we were informed
of \cite{ppp} which does the same for the v section. Although the table
does not tell us which of these knots is fibred (and now does not need to, in
light of this section and Dunfield's list), we find in it eight of our 
unknown 1-cusped 3-manifolds including the last three to be dealt with.
The descriptions given of these three knots are: v3093 is 16n245346 in
Knotscape (if it had been an alternating knot then our work would have been
in vain because we would have been able to conclude that it was fibred just
from the Alexander polynomial). Then v2869 and v3541 are given in terms of
a (non-alternating) Dowker-Thistlethwaite code with 18 and 21 crossings
respectively. Although these may not be the minimal crossing numbers, they
must come pretty close because Knotscape tells us they are not in its census
which goes up to 16 crossings. Also we now know the topological type of
their fibres, because as knots in $S^3$ their fibres will have one boundary
component and genus half the degree of the Alexander polynomial.

In conclusion we have:
\begin{prop}
The proportion of fibred 3-manifolds in the (orientable) cusped census is
exactly 4199/4815=0.87206645898...
\end{prop} 

\section{Virtually fibred cusped 3-manifolds} 
As we now know all fibred 3-manifolds in the cusped census, we turn to how
we can find non-fibred virtually fibred examples. The crucial point is that
a non-fibred hyperbolic 3-manifold that is commensurable
with a fibred hyperbolic 3-manifold is itself virtually fibred, by considering
the common finite cover, so that the property of being virtually fibred is
constant on commensurability classes. Therefore we ought in principle to be
able to use our fibred 3-manifolds to obtain non-fibred commensurable
examples $M$. The first case that comes to mind is when
$\pi_1M$ is arithmetic, which in the cusped case
means that it has integral traces and the
invariant trace field is an imaginary quadratic number field. Here two
arithmetic fundamental groups will be commensurable if they have the same
invariant trace field, so on finding a fibred example we have that all
arithmetic hyperbolic cusped 3-manifolds with this imaginary quadratic
number field will be virtually fibred.

However recently the paper \cite{ghh} gives an algorithm that determines the
commensurator of any non-arithmetic cusped hyperbolic 3-manifold and it is
then applied to find commensurability classes for the 3-manifolds in the cusped
census, as well as for hyperbolic knots and links for up to twelve crossings.
Therefore it is worth looking at the 616 non-fibred census 3-manifolds to
see if any are in the same commensurability class as a fibred 3-manifold,
given that we now can recognise all fibred 3-manifolds in the cusped census.
Doing this gives us 86 non-fibred virtually fibred cusped hyperbolic
3-manifolds as listed in Table 2 (a few of which would have been known
before, see for instance \cite{cd} and \cite{hmw}). Most of the fibred
3-manifolds certifying that these examples are virtually fibred have more
than one cusp; moreover the four non-fibred 3-manifolds with 2 cusps
(v2943, v3379, v3384, v3396) all appear thus we can say that any hyperbolic
3-manifold in the census with more than one cusp is virtually fibred.

We can further add to this table because the data we are using includes
commensurability classes of knots and links in $S^3$. However, rather than
just looking for fibred knots and links, we use the recent result \cite{wal}
that all 2-bridge knots and links are virtually fibred. We can identify
2-bridge knots and links in the tables by their Conway notation.
This gives us another 51 examples to add to our table. Most of these are
themselves non-fibred 2-bridge knots or next to one in the census, 
although a few are shown virtually fibred by being commensurable 
with a 2-bridge knot
that is not in the cusped census. We have also two links
not from the census that make an appearance: there is the fibred 2-bridge
link 8a31 (or $8^2_4$ in the tables) with Conway notation 323 and the
non-fibred 2-bridge link 10a171 with Conway notation 262 (in fact the
2-cusped 3-manifolds v2943 and v3379 mentioned above are also 2-bridge
links identifiable as 7a11 or $7^2_3$ or 232 and 8a24 or $8^2_6$ or 242
respectively).

One amusing consequence of the ubiquity of 2-bridge knots amongst those
with low crossing number is that just by striking out from the tables of
knots with nine crossings or less the 2-bridge knots and the knots with
monic Alexander polynomial (which for these crossing numbers will be
fibred), we see that the only ones left that are not known to be
virtually fibred are the ten knots $8_{15}$ (8a2), $9_{16}$ (9a25), 
$9_{25}$ (9a4), $9_{35}$ (9a40), $9_{37}$ (9a18), 
$9_{38}$ (9a30), $9_{39}$ (9a32), $9_{41}$ (9a29), $9_{46}$ (9n5) 
and $9_{49}$ (9n8). There may be a few more cusped
3-manifolds in the census that could be added to this table by having full
knowledge of which knots and links up to twelve crossings are fibred, but
certainly some non-fibred 3-manifolds are listed alone in their
commensurability class so this process would not finish the job off.
However we have pushed the number of virtually fibred 3-manifolds in 
the cusped census up to 4336 which is a fraction over 90\%. 

\section{Closed fibred hyperbolic 3-manifolds}
In the Hodgson-Weeks census \cite{hw}
of closed hyperbolic 3-manifolds, consisting of
just under 11,000 examples (the number given is 11,031 but there are a few
duplications), nearly all have finite first homology: only 127 have first Betti
number 1 and above that there is but one 3-manifold with first Betti number 2.
Thus only these few special closed 3-manifolds have a chance of being fibred,
but in fact there is a reason why it is likely to be a good chance. All
3-manifolds in the closed census are obtained by Dehn surgery on 1-cusped
3-manifolds from the cusped census and this process either preserves the
first Betti number or reduces it by one. Therefore the closed 3-manifolds $M$
with $\beta_1(M)=1$ come from 1-cusped 3-manifolds $M'$ with $\beta_1(M')=1$
or 2. But there are only 3 examples of the latter and moreover we now know
that the vast majority of 3-manifolds $M'$ in the cusped
census are fibred. If so and if $\beta_1(M')=1$ then we have mentioned 
in Section 2 that $M$ must be fibred too. 

In addition the one closed 3-manifold $M$ with $\beta_1(M)=2$ happens to be
v1539(5,1), so it is irreducible and therefore
Section 2 tells us it is fibred, as well as v1539(-5,1) which also appears in
the census. Otherwise we work through the closed 3-manifolds $M$ with
$\beta_1(M)=1$, seeing if they are surgery on a 1-cusped 3-manifold $M'$
that is listed as fibred but which is not one of the 
three special cases with $\beta_1(M')=2$. In this way we find 80 further
closed fibred 3-manifolds in the census which is a big proportion
of those with positive first Betti number. The results are listed in Table 3.

As for the remaining 46 closed 3-manifolds $M$ with
$\beta_1(M)>0$ in the census, we calculate the Alexander polynomial of the
given fundamental group presentation which proves
that all but five are not fibred. As we have $\beta_1(M)=1$, we can do this
in exactly the same way as we did for 1-cusped 3-manifolds, and indeed 
it is still invariant under $t\mapsto t^{-1}$. 
Moreover it is again the case that if $M$ is fibred over the
circle then $\Delta_M$ must be monic, and here the degree of $\Delta_M$ must
be twice the genus of the fibre: we can see this from (2) by noting that we 
need to add a relation for the closed surface, but this results in an extra
row of zeros on application of the free differential calculus.

Our fundamental groups
are usually 2 generator, 2 relator with a few 3 generator, 3 relator examples
but we can use short cuts that might avoid calculating the whole Alexander
polynomial. If we have $\pi_1M=\langle g,x|r_1,r_2\rangle$, which we always
assume is in standard form with respect to $x$, then $\partial r_i/\partial x
=0$, thus the Alexander polynomial is 
the highest common factor of the two polynomials
$\partial r_i/\partial g$. But as we know $M$ is hyperbolic, if it is fibred
then this must be by a surface of genus at least two, so the Alexander
polynomial must be monic of even degree at least four. We thus calculate
only one polynomial corresponding to the nicest looking relation and if this
does not have such a factor then we are done. It turns out, 
as seen in Table 4, that in all but three of the cases the polynomial
obtained was quartic, non-monic and not a scalar multiple of a 
monic quartic
polynomial, so these 3-manifolds are not fibred. The three exceptions were
that with v2018(-4,1) a quintic was obtained which factors as 
$(t+1)(t^2+1)(2t^2-3t+2)$ so this is non-fibred, indeed the other relation 
gives $(t^2+t+1)(t^2+1)(2t^2-3t+2)$
so the last two factors are the Alexander polynomial. This 3-manifold
will feature again in Section 7 where we will find that it is virtually
fibred. The next 
exception that needs to be checked is v2238(-5,1), but here
a quintic is obtained that factors into irreducibles as $(t+1)(2t^4-t^3-t+2)$
so this is fine. The only other problem is v3183(-3,2) which yields $2(t^4+1)$
so we worry that $t^4+1$ might be the Alexander polynomial, but 
looking at the other relation we see this cannot be the case.

As for the three 3 generator cases, we similarly take 2 relations and
calculate the relevant $2\times 2$ determinant; these are all quartic and
present no problems. We treat those closed 3-manifolds which come from the
three special 1-cusped 3-manifolds s789, v1539, v3209 separately. For the
2 generator group $\pi_1(\mbox{v}1539)$ 
we have already stated in Section 2 that
$(Ab,B^3a^5B^2)$ is a basis for the cusp, so taking the relation
$(Ab)^p(B^3a^5B^2)^q$ from v1539($p,q$) 
and substituting $a=bx$ so that it is in
standard form with respect to $b$ gives us the polynomial
\[q t^4+qt^3+(q-p)t^2+qt+q\]
whereas the original relation gives 0, so this is the Alexander polynomial
(except for $(p,q)=(5,1)$ where $\beta_1(M)=2$)
and $q\neq 0,1$ implies that the 3-manifold
is not fibred.
We now have built up the complete picture for these hyperbolic 3-manifolds as
we saw in Section 2 that v1539$(p,1)$ is fibred (and it is clear that
v1539$(1,0)$ has cyclic fundamental group so is not hyperbolic);
in particular v1539(5,2)
that appears in Table 4 is non-fibred. Similarly for s789 we have
$(abc^2,a^3cbcA^3C)$ as a basis for the cusp and we take this Dehn filling
relation for s789$(p,q)$ along with either one of the two original
relations (they result in the same polynomials).  
We put $c=Ax$ and $b=ya$ to get two relations in standard form
with respect to $a$ and this yields the Alexander polynomial
\[qt^4-qt^3+(p+q)t^2-qt+q\]
so once again it is not fibred if $q\neq 0$ or 1 
(with $\pi_1\mbox{s}789(1,0)=\mathbb Z$ again),
sorting out s789(-5,2). Finally
we do this for v3209, with basis $(aCbc^2,aCacAcAC)$ and either one of
the original relations, setting $a=Cx$ so that we are
in standard form with respect to $c$. For v3209$(p,q)$ we have
$\pi_1\mbox{v}3209(1,0)=\mathbb Z$ and Alexander polynomial
\[qt^4-2qt^3+(p+2q)t^2-2qt+q\]
which reveals nine closed 3-manifolds in Table 4 as not fibred when
$q>1$. 

We guess that s789$(p,1)$ and
v3209$(p,1)$ are all fibred; not only would this fit into the same pattern
as v1539 
but we have already seen in Table 3 that s789$(p,1)$ for $p=\pm 5$
and v3209$(p,1)$ for $p=\pm 3$ are fibred as they have alternative descriptions
as Dehn fillings on 3-manifolds $M$ with $\beta_1(M)=1$. 
We can say that if so, they must have fibres of genus two. 

However this still leaves in the census five 3-manifolds 
v3209$(p,1)$ for $p=\pm 4,\pm5,6$
whose status is unknown. In the hope of finishing this off, it is worth
looking for cyclic covers which we can show are fibred, just as we did with
the remaining 1-cusped 3-manifolds in Section 4. Happily this works for
all five thus the fibred status of every 3-manifold in the closed census is
known: 87 are fibred, 41 are non-fibred with $\beta_1(M)=1$ and the rest are
non-fibred with $\beta_1(M)=0$. We summarise the details so as to allow the 
claims to be checked. All five cases are very similar. We put $a=xC$ in
our presentation and then we have fundamental group $\langle a,c,x\rangle$ in
standard form with respect to $c$. We know the fibre would be a genus 2 surface
so we are after a 5 generator presentation. In each case the cyclic covers
of degree 2 and 3 have too few generators (at least on rewriting) but Magma
tells us that the cyclic cover of degree 4 yields a 5 generator presentation
of the form $\langle g_1,g_2,g_3,g_4,t\rangle$ for
\begin{eqnarray*}
(g_1,g_2,g_3,g_4,t)&=&(x,cxC,Cbc,c^2xC^2,c^4)\\
&&(cbC,cxC,Cbc,c^2xC^2,c^4)\\
&&(x,cxC,Cbc,Cxc,c^4)
\end{eqnarray*}
where the first option is for $p=4,\pm 5$, the second for $p=-4$ and the third
for $p=6$. As $t=c^4$ has infinite order but all $g_i$ have finite
order in the homology of $M$,
we know the presentation obtained in each case will be in
standard form with respect to $t$. What is most promising is that we always
find the first relation given has no appearance of $t$ at all (but $t$ does
appear in the others). Indeed in all but $p=-4$ this relation is of length
8 with each $g_i^{\pm 1}$ appearing once, which is a relation defining the
closed surface of genus 2. For $p=-4$ it is of length 12 but as a consequence
of showing the 3-manifold is fibred, this relation has to define the genus 2
closed surface group as well.

We then proceed just as in Section 4 by looking at the subwords from $t$ to
$T$, or from $T$ to $t$ (we did in fact do both).
In all but $p=-4$ we are given more than 5
relations so we are looking for generating sets for the free group on
$g_1,g_2,g_3,g_4$ rather than a free basis, but we always proceed by
taking our $n$ subwords (where $n$ can be 4, 5 or 6) and using the shorter
subwords to knock letters off the longer subwords until we have each
generator $g_i$. We do this by hand: for $p=\pm 4$ the relations are in
simple form with respect to $t$. For $p=5$ the fourth and sixth of the seven
relations have two appearances of $t$ (whereas
the first relation has none and the rest have one). They are of the form
$tw_1Tw_2tu_1Tu_2$ and $v_1tv_2TW_2tW_1T$ for $u_j,v_j,w_j$ words in the
$g_i$ so we can concatenate them to obtain a relation in simple form which
we now use. For $p=-5$ we have six relations with the third, fifth and sixth
in this double form but each pair of these three can be concatenated as above
to obtain five relations in simple form. Then for $p=6$ we are given seven
relations with the last three simple. We put together the second and fifth
to obtain $tsT$, where $s=Cxc$, which we can now insert into the three
relations in double form, resulting in enough relations in simple form to
obtain all the generators.

Finally to show the original 3-manifolds are fibred, we look at the homology
of the degree 4 covers. These are listed below and all have first Betti
number 1 so we are done. 

\begin{table}[h]
\centering
\begin{tabular}{lcl}
3-manifold&\qquad&Homology of cover\\
v3209(4,1)&\qquad&$\z_2+\z_4+\z_4+\z_{24}+\z$\\
v3209(-4,1)&\qquad&$\z_2+\z_4+\z_4+\z_8+\z$\\
v3209(5,1)&\qquad&$\z_5+\z_5+\z_{65}+\z$\\
v3209(-5,1)&\qquad&$\z_5+\z_5+\z_{15}+\z$\\
v3209(6,1)&\qquad&$\z_2+\z_6+\z_6+\z_{42}+\z$\\
\end{tabular}
\end{table}

Thus we now know all the fibred 3-manifolds in the closed
census. We have seen that if $M'$ is a 1-cusped fibred 3-manifold
with $\beta_1(M')=1$ and we Dehn fill along its longitude to create $M$ then
$M$ is fibred. We might expect that if instead $M'$ is non-fibred then $M$ is
not but this is unlikely to be true in full generality. For instance 
let us take the 1-cusped 3-manifold m137 (an interesting example as it has
a quadratic imaginary invariant trace field but is the first in the cusped
census not to have integral traces). It is not fibred (indeed is not known
to be virtually fibred) and is a knot in an integral homology sphere. 
We find from SnapPea a fundamental group presentation and basis for the
cusp, whereupon it is easily seen that the group $\z$ is obtained on Dehn
filling of the longitude thus (assuming Poincar\'e) $M=S^2\times S^1$ and so is
fibred. (Another 3-manifold $M'$ in the census with $\beta_1(M')=1$ where
$\mathbb Z$ is obtained on Dehn filling is the non-fibred s783, 
as well as the three 1-cusped examples with $\beta_1(M')=2$.)
However if $M'$ is the exterior of a non-trivial knot in $S^3$ then
Gabai shows in \cite{gab} that $\pi_1M\neq\z$. He goes on to prove that for
knots $M'$ is fibred if and only if $M$ is, in which case the fibres have the
same genus. Although this seems useful, and certainly we have included in Table
3 the genus of the fibre of those closed 3-manifolds $M$ where the given
$M'$ is a knot exterior in $S^3$, 
there was only one case where this would have proved $M$ is
non-fibred: s862 is the non-fibred knot $8_4$ so s862(7,1) in Table 4 is not
fibred. In trying to generalise Gabai's result,
a conjecture of Boileau (Problem 1.80 (C) in the
Kirby problem list \cite{kir}) states that if $K$ is a null-homotopic
knot in a closed orientable irreducible 3-manifold $M$ then a non-trivial
Dehn surgery on $M-K$ produces a fibred 3-manifold if and only if $M-K$ is
fibred and it is the longitudinal surgery. Here the trivial surgery is just
filling in $K$ to obtain $M$ thus destroying the meridian, 
and a null-homotopic knot can be detected because the longitude then 
becomes trivial. A fair variant on this question might be:
if $M'$ is a 1-cusped hyperbolic 3-manifold with $\beta_1(M')=1$
where the longitudinal surgery produces a closed 
fibred 3-manifold $M$ that is hyperbolic then is
$M'$ fibred? This is true for all examples we have considered.

\section{Virtually fibred closed 3-manifolds}

We will now use our data to find non-fibred virtually fibred closed hyperbolic
3-manifolds. There seem to be even less examples of these than in the cusped
case: until this point the only known ones in the literature consisted of 
the original idea due to Thurston of the union of two
twisted $I$-bundles over a non-orientable surface, which have a fibred
double cover, and the pair of non-Haken examples in \cite{rei} (one of
which is the unique double cover of the other). However, just as in the
cusped case, we merely need to find non-fibred hyperbolic 3-manifolds that 
are commensurable with fibred hyperbolic 3-manifolds. 
In particular any 3-manifold $M$ in the closed
census which is commensurable with something in Table 3, but which is not in
Table 3 itself, is a non-fibred virtually fibred example. 
We certainly do not have a full enumeration of the commensurability classes
as in the cusped case, so we turn to the theory of
arithmetic Kleinian groups: that is if we have arithmetic hyperbolic
3-manifolds $M_1,M_2$ then they are commensurable if and only if their
invariant trace fields and invariant quaternion algebras are isomorphic. In
the closed arithmetic case we are guaranteed more invariant trace fields
than just the imaginary quadratic ones: in fact the fields that occur
are precisely those with exactly one conjugate pair of complex embeddings. In
order to determine this we utilise the program Snap \cite{sna} (see
\cite{cghn} for a description) and look for the file 
\texttt{snap\_data/closed.fields} which lists (in order of volume) all
closed 3-manifolds in the closed census for which the invariant trace 
field and invariant quaternion algebra could be found. It is
known that $M=\mathbb H^3/\Gamma$ is arithmetic if and only if the invariant
trace field $k\Gamma$
has exactly one complex place, the invariant quaternion algebra $A\Gamma$
is ramified at every real place and $\Gamma$ has integer traces. Thus if
$M$ is a fibred 3-manifold from Table 3 appearing in this list we next look
at the file \texttt{snap\_data/closed\_census\_algebras} which gives (listed
in order of trace field) 3-manifolds grouped together by invariant trace field,
quaternion algebra, and whether or not they are arithmetic. Hence if $M$ is
arithmetic then
all 3-manifolds appearing together in the same grouping as $M$ are 
commensurable with $M$, and so virtually fibred.

The results are listed in Table 5. In particular we find that the 
Weeks 3-manifold
m003(-3,1), conjectured to be the smallest volume closed hyperbolic 3-manifold
and known \cite{chin} to be the smallest volume arithmetic 3-manifold, is
virtually fibred as it is commensurable with m289(7,1). The third entry
m007(3,1) in the closed census is one of the two non-Haken virtually fibred
closed 3-manifolds in \cite{rei} and is called Vol(3) as it is
the conjectured third smallest closed hyperbolic 3-manifold. This is known
to be arithmetic (see \cite{jr}) so we can add it and the other 3-manifolds
that Snap lists in its commensurability class to Table 5. Work of Dunfield
\cite{dunhak} determines that out of the 246 3-manifolds in the closed census
with volume less than 3, exactly 15 are Haken. Only one from that list
appears here (this is m140(4,1) with volume 2.6667)
so all other 3-manifolds in Table 4 with volume less than 3 are
non-Haken virtually fibred hyperbolic
examples. For other specific examples of Haken non-fibred
virtually fibred closed hyperbolic
3-manifolds, one can use Theorem 2 in \cite{rei}
which shows that the $3k$-fold cyclic branched cover $M_{3k}$ of the
figure eight knot is a double twisted $I$-bundle with $\beta_1(M_{3k})=0$. 
However we also have, as promised, a closed non-fibred virtually fibred 
3-manifold in the form of v2018(-4,1) with positive Betti number. 
Incidentally it can be checked that this
3-manifold is genuinely a new example and not a union of two twisted 
$I$-bundles because if so it would have a fibred double cover, but
all its three index 2 subgroups have first Betti number 1.
We claim that this is the first known example of its kind:
for instance in \cite{bw} it is shown that for every $n>0$
there exist non-fibred closed hyperbolic 3-manifolds $M_n$ with 
$\beta_1(M_n)=n$ but it is not known if they are virtually fibred.  

We end up with 129 non-fibred virtually fibred 3-manifolds from the closed
census. One might say that this is only a small proportion of the whole
census, but of course our method only gives rise to arithmetic examples because
$(k\Gamma,A\Gamma)$ is not a complete commensurability invariant in 
the non-arithmetic case. Another point is that
all the examples of virtually fibred 3-manifolds we have given are 
commensurable with fibred 3-manifolds that necessarily must appear
in the census, whereas as the volume grows and we have more and more
3-manifolds one would expect to have to look further for commensurable fibred
3-manifolds. This could explain why we do better with the
3-manifolds of smallest volume: of the first 51 census 3-manifolds 
(which goes up to volume twice that of the regular ideal tetrahedron), 34
are arithmetic, with 15 of these now known to be virtually fibred.

\section{Co-rank of the census 3-manifolds}
The co-rank $c(G)$
of a finitely generated group $G$ is 
the maximum $n$ for which there is a homomorphism from $G$ onto the free
group $F_n$ of rank $n$. Clearly $\beta_1(G)\geq c(G)$ and 
$\beta_1(G)\geq 1$ implies $c(G)\geq 1$. This quantity is of algebraic
interest and we can think of the property $c(G)>1$ as giving rise to one of
the several notions of ``largeness'' of a group; see for instance
\cite{butjpaa}. But if
$G=\pi_1M$ for $M$ a compact orientable 3-manifold (for which we write $c(M)$)
then we have a geometric interpretation which allows us to think of it as a
measure of ``largeness'' of a 3-manifold: this is because $c(M)$
is the maximal number of disjointly and properly embedded orientable connected
surfaces $S_i$ for which $M\backslash \cup S_i$ is 
connected (and in this context is also called the cut number of $M$).  
We can ask about the co-rank of 3-manifolds in the closed or cusped
census: this can quickly
be determined for every single one, and it turns out that we do not have
any examples of ``large'' 3-manifolds here. As pointed out in \cite{htrs},
there is a (computationally very inefficient) procedure to determine if
a finitely presented group surjects onto $F_n$, but it will not prove the
non-existence of such a surjection. However, in this setting we have available
properties of 3-manifold groups to help us.   
\begin{theorem}
If $M$ is a 3-manifold appearing in the closed census then $c(M)=0$ if
$\beta_1(M)=0$ and otherwise $c(M)=1$. If $M$ is a 3-manifold
appearing in the cusped census then $c(M)=1$.
\end{theorem}
\begin{proof} We only need to do anything when $\beta_1(M)>1$. However if
so and if $M$ is fibred then $\beta_1(M)>c(M)$. This is Theorem 4.2 
in \cite{but} but here is a 
variation on that proof. If $\beta_1(M)=c(M)=n$ with 
$\theta:\pi_1M\rightarrow F_n$ a surjective homomorphism 
then any homomorphism from $\pi_1M$ to $\mathbb Z$ factors
through $\theta$. If $M$ is fibred then we 
have our finitely generated kernel $K$ of our relevant surjective homomorphism
in $\pi_1M$ which is normal and of infinite index, so $\theta(K)$ 
must be be the
same in $F_n$. But non-abelian free groups do not have finitely generated
normal subgroups of infinite index except for the trivial group.

Thus this sorts out v1539(5,1), the only closed 3-manifold with Betti number
2. It also sorts out all cusped 3-manifolds $M$ 
(which must have $\beta_1(M)\geq 1$) except for the four non-fibred examples in
Section 4 Part 2 with $\beta_1(M)=2$ and the three fibred examples in the
census with $\beta_1(M)=3$. For these seven, we have to eliminate the
possibility that $c(M)=2$.

Firstly v2943 and v3379 are 2 generator, so we cannot have $\pi_1M$ 
surjecting onto $F_2$ unless $\pi_1M=F_2$ which is not true. The given
presentation for $\pi_1$(v3384) is 
\[\langle a,b,c| ab^2ab^2aCb^2ab^2abcb, \,\,aCAc\rangle.\]
The second relation means that our surjection $\theta$ onto $F_2$ would 
have to send
$a$ and $c$ onto powers of the same element $v\in F_2$ because that is the
only way elements can commute in a non-abelian free group. So $u=\theta(b)$ and
$v$ must generate $F_2$, hence be a free basis, but this is not possible by
looking at the image of the first relation which would always
give a non-trivial relation between $u$ and $v$.

This argument also works for the three 3-manifolds s776, v3227, v3383 with
$\beta_1(M)=3$: we know $c(M)=3$ is not possible and to eliminate $c(M)=2$
we use the second relations given in each case. Respectively they are
$aCAc$, $bCBc$, both of which work in exactly the same way above, and
$aCb^2AcB^2$, which by setting firstly $a=cx$ and then $c=b^2Y$ becomes
$b^2YxyXB^2$, so we now just use the pair of generators $x,y$.

This leaves only
\[\pi_1\mbox{(v3396)}=\langle a,b,c|aBca^2bC,\,\,a^2cba^2CAB\rangle\]
with abelianisation $\mathbb Z_3+\mathbb Z+\mathbb Z$. We suppose 
$\theta:\pi_1$(v3396)$\,\rightarrow F_2$ is onto and to finish we
derive three quick contradictions.
Both groups have three subgroups of index 2, which in the case of $F_2$
are all copies $H_i$ of $F_3$. As each $\theta^{-1}(H_i)$ is distinct and has
index 2, these must be the three index 2 subgroups $K_i$ of $\pi_1$(v3396)
so $c(K_i)\geq 3$, which implies that $\beta_1(K_i)\geq 3$ and $K_i$ will
need at least four generators. Two subgroups pass those tests but the third
is $\langle a,cb^{-1},b^2\rangle$ and has abelianisation $\mathbb Z_{24}+
\mathbb Z+\mathbb Z$ so it fails on both counts. Or we could try the lazy
approach: by considering $\theta^{-1}(H)$ for $H$ finite index in 
$F_2$ as before
we have that $\pi_1$(v3396) must have as many subgroups of index $n$ as $F_2$
does, so we ask the computer. The numbers we get 
from index 2 onwards are 3,15,32,64 for $\pi_1$(v3396) whereas for $F_2$ 
they are 3,7,26,97 so we have already been overtaken at index 5. In fact
this is actually the number of subgroups up to conjugacy but our point
still holds.
\end{proof}

\section*{Appendix: Guide to Tables}
\markboth{APPENDIX: GUIDE TO TABLES}{APPENDIX: GUIDE TO TABLES}
{\bf Table 1}: Alexander polynomials of unknown cusped census 3-manifolds\\
{\bf Table 2}: Cusped virtually fibred non-fibred census 3-manifolds\\
{\bf Table 3}: Closed fibred census 3-manifolds\\
{\bf Table 4}: Closed non-fibred census 3-manifolds with infinite homology\\
{\bf Table 5}: Closed virtually fibred non-fibred census 3-manifolds\\
\hfill\\
{\bf Notes on Tables}\\
\hfill\\
{\bf Table 1}: This lists in the column ``Name'' the 165 cusped 3-manifolds $M$
with $\beta_1(M)=1$ which are unknown in Dunfield's list\\
\texttt{http://www.its.caltech.edu/\textasciitilde dunfield/snappea/tables/}\\
\texttt{mflds\_which\_fiber}
of fibred and non-fibred cusped 3-manifolds. For each one, we take the 
presentation for its fundamental group (as given in\\ 
\texttt{virtual\_haken\_data/manifolds/cusped.gap}
available at\\
\texttt{http://www.its.caltech.edu/\textasciitilde
dunfield/virtual\_haken/})\\ 
which is always
(with the exception of v3036 which is marked by *2 gen*) generated by $a,b,c$
and with two relations. The ``Standard column'' indicates the substitutions we
must make, in order, to put the presentation into standard form with respect
to a generator (meaning that the generator has zero exponent sum in both
relations); this generator is then given at the end. Then the column ``Poly''
gives the Alexander polynomial which is written in a compact form. If a single
number $n$ is given without brackets then the presentation obtained was in
simple form, as described in Section 4 Part 3, so that the Alexander polynomial
must be of the form $nt+m+nt^{-1}$. Here $n$ can be obtained quickly and we
do not need to calculate $m$, unless $n$ is zero in which case we do and we
write $0=[m]$. The brackets notation that we use in general is because
the Alexander polynomial is equal, up to units, when $t$ is substituted for
$t^{-1}$ and it is non-zero when evaluated at 1. Thus it is either of even
degree and in the form
\[a_kt^k+\ldots +a_0+\ldots +a_k^{-k},\mbox{ written }[a_k,\ldots ,a_0]\]
or of odd degree in the form
\[a_kt^k+\ldots +a_1t+a_1+\ldots +a_k^{-(k-1)},\mbox{ written }
(a_k,\ldots ,a_1).\]
The six 3-manifolds that have monic Alexander polynomial are printed in bold,
as is the leading coefficient. They are all fibred except v2530.\\ 
\hfill\\
{\bf Table 2}: Here we list under ``Name'' the non-fibred virtually fibred
cusped census 3-manifolds that we found (we know they are non-fibred by
Dunfield's list and the results of Section 4) using the file of cusped
commensurability classes that makes up the data resulting from \cite{ghh} 
(supplied to us by the authors, for which we thank them). In the column
``Name of fibred'' we list the fibred 3-manifolds
with which the listed 3-manifolds are commensurable, 
thus showing that they are virtually fibred. The column before this is headed 
``Ratio'' and is the
ratio of the volume of the virtually fibred 3-manifold(s) to that
of the corresponding group of fibred 3-manifolds. The 3-manifolds with
{\bf 2} or {\bf 3} as a superscript have that number of cusps whereas the
rest all have one cusp. As mentioned in Section 5, we also use 2-bridge
knots and links. Here several notations are in use, so we give its name as a
census 3-manifold (if it is one) as obtained from \cite{cdw} and \cite{ppp},
then the Knotscape name (crossing number, a (or n) for (non-)alternating
and the reference number) then the ordering in the knot tables started by
Alexander and Briggs, and extended by Rolfsen and Bailey using work of
Conway. This only applies for knots with ten or less crossings and links
of nine or less. Then we give the Conway notation, needed to confirm it is
2-bridge, in which case this is just a string of integers (written together,
with two digit numbers denoted [10] etc).

In order to move between these different notations,
the file has commensurability classes of knots and links up to twelve
crossings given under the Knotscape name, which it lists as equal to the
relevant cusped census 3-manifold if appropriate. For knots of 10 crossings
or less we can use the file in Knotscape that converts between its notation
and the Rolfsen-Bailey tables, then look up the Conway notation in \cite{rolf}.
For 11 crossing alternating knots, the original enumeration is due to Little
but it was then taken up by Conway. We found\\
\texttt{http://www.indiana.edu/\textasciitilde knotinfo/}\\
which converts from Knotscape to Conway notation. To check this, we then
have\\
\texttt{http://www.scoriton.demon.co.uk/knots.html}\\
which allows us to go from Conway notation to braid notation (this table is
in order of Little's notation so we confirm it with Conway in \cite{con})
which we can then enter into Knotscape and ask it to identify the knot, thus
taking us back.

There was one census knot each for 12 and 13 crossings that featured; by
getting Knotscape to draw them it was immediately seen that they were
both twist knots. For the two links, we used \cite{ad} to go between
Thistlethwaite's notation as given in the file and the Rolfsen-Bailey tables
by recognising volumes in one case, whereas for the ten crossing link we
recognised it as a 2-bridge link from the picture in\\
\texttt{http://www.math.toronto.edu/\textasciitilde drorbn/KAtlas/Links/}\\
 
Finally non-fibred arithmetic 3-manifolds are confirmed virtually fibred by
the symbol $\mathcal An$ in the ``Name of fibred'' column, where $n$ can
be 1,2,3 or 7 which refers to the imaginary quadratic number field which is
its invariant trace field. As we know of arithmetic fibred cusped 3-manifolds
with each of these invariant trace fields, 
they will be commensurable with those listed under ``Name''.\\
\hfill\\ 
{\bf Table 3}: This lists all 
closed 3-manifolds in the census which are fibred, as shown 
in Section 6. There are 87 entries listed in order of volume, which is
given in the first column as it can be time consuming to find a 3-manifold
by hand on name alone. To aid this, the volume is given to 4 decimal places,
which should be enough to find the right part of the census, and is always 
rounded down to avoid having to look back. The \y \ symbol indicates a volume
which is the same as the preceding volume to the accuracy given in the census.
Next we give the name of the 3-manifold as listed in the census, which we
take to be\\
\texttt{ftp://www.geometrygames.org/priv/weeks/SnapPea/SnapPeaCensus/}\\
\texttt{ClosedCensus/ClosedCensusInvariants.txt}\\

The column ``$\mathbb Z$'' refers to those 3-manifolds whose homology is
$\mathbb Z$ and a dot indicates this. If the associated cusped 3-manifold 
is a knot in $S^3$ (as given by \cite{cdw} and \cite{ppp})
then the corresponding closed 3-manifold is then surgery along a longitude 
so its fibre will have the same genus as the knot,
in which case we put this number in the column instead. As shown in Section
6, the genus of the fibre of any of the 3-manifolds in this table can be
calculated from the fundamental group presentation
if required. The $\beta 2$ indicates the one 3-manifold with
homology $\mathbb Z+\mathbb Z$. The ``neg'' column 
marks with - those 3-manifolds
that are listed in the census as having negatively oriented tetrahedra
present. The program SnapPea has alternative descriptions for some 3-manifolds
which might not involve negative orientations. In the ``Alternative'' 
column we
have included such a description in one case, as well as alternative
descriptions known to us for 3-manifolds obtained by $(p,1)$ surgery on the
3-manifolds s789 and v3209 as this is required to prove they are fibred.\\
\hfill\\
{\bf Table 4}: This lists the remaining 41 closed 3-manifolds 
in the census with
infinite homology, along with evidence to show that they are non-fibred. They
are given by volume and name, then in the ``Standard'' column we give the 
substitutions we used to put their fundamental groups in standard form,  
followed by the relevant generator, starting from the presentations given in
\texttt{virtual\_haken\_data/manifolds/final.gap}\\
at 
\texttt{http://www.its.caltech.edu/\textasciitilde
dunfield/virtual\_haken/}\\ 
We note in this
column that s862 is the (non-fibred) knot $8_4$. In the ``Poly'' column we give
the polynomial obtained from the first relation, using the same notation for
polynomials as in Table 1 (so the Alexander polynomial
is a factor of this but we have not confirmed that they are equal). From
Section 6 this polynomial immediately tells us that the 3-manifold is 
non-fibred except for the three indicated in bold for which we refer back
to that section.

For a 3-manifold that is $(p,q)$ surgery
on s789, v1539 or v3209, we show in Section 6 that
$q\neq 1$ implies it is non-fibred and so we mark these with \texttt x.\\
\hfill\\
{\bf Table 5}: This lists the closed virtually fibred 3-manifolds found in
Section 7; they are all arithmetic. Also they all have finite homology (hence
are non-fibred) with one exception, marked by the suffix $\beta 1$ and
printed in bold. Again we list
volume, name (at 2.5689 we list m130(-3,1) with ? because it is given as
m130(1,3) in the original census but the former in all other sources) and the
column ``neg'' marks those 3-manifolds with negatively oriented tetrahedra (at
this point we did not have access to possible alternative descriptions).

As in Table 2 for the cusped case,
in the column ``Name of fibred'' we give the fibred 3-manifolds from Table 3
with which the listed 3-manifolds (put together in a group if they are
commensurable and have the same volume) are commensurable, thus showing that
they are virtually fibred. There is one commensurability class that is
proved virtually fibred by using Vol(3) in \cite{rei} which is in the
census as m007(3,1). We put a zero superscript on this to remind ourselves
it has zero first Betti number. We then have in ``Ratio'' the
ratio of the volume of the virtually fibred 3-manifolds in each group to that
of the corresponding group of fibred 3-manifolds (it happens that
the latter always have the same volume within a group). They are given as
fractions with small coefficients; although this is likely to be correct, it
could be argued that unlike in the cusped case where we are able to
use the index of the 3-manifold in its commensurator
we have only confirmed it to the number of decimal
places available. This does not concern us here because the aim is to allow
quick access to the volumes of those 3-manifolds in the right hand column for
ease of reference.\\
\hfill\\

\newpage
\markright{APPENDIX: TABLES}

\begin{table}
\caption{Alexander polynomials of unknown cusped census 3-manifolds}
\centering

\begin{tabular}{|lrr|lrr|}
\hline
Name&Standard&Poly&Name&Standard&Poly\\
\hline
m306&$c=xA;a$&-3&m307&$b=Ax;a$&[3,2]\\
m372&$a$&-2&m373&$a$&-2\\
m410&$a$&[-2]&&&\\
\hline
s386&$a$&-2&s387&$c$&-2\\
s426&$b=Ax;a$&[-4,2]&s427&$b=Ax,c=Ay;a$&[-4,-2]\\
s435&$a$&-3&s436&$c=Bx;b$&3\\
s486&$b=xA;a$&[-5,4]&s487&$c=xA;a$&[-5,-4]\\
s491&$c=yA,b=Ax;a$&[-5,-6]&s492&$b=xA;a$&[5,-6]\\
{\bf s594}&$b=Ax;a$&({\bf 1},3)&s626&$b=YX^3,c=x^2y;x$&[3,-3,2]\\
s673&$a=x^2y,b=YX^3;x$&[3,-3,4]&{\bf s704}&$c=Z^3y,a=Yz^2;z$&[{\bf  1},-2,3]\\
s707&$b$&-3&s708&$b=Ax;a$&[3,0]\\
s732&$c$&[-2,-11]&s733&$a$&2\\
s773&$c=Ax;a$&(2,2)&s779&$c=Ax;a$&(2,0)\\
s784&$b$&3&s788&$c$&-3\\
s818&$a$&-2&s819&$a$&-2\\
s837&$c$&-3&s838&$a$&3\\
s878&$b$&2&s$879$&$c=xA;a$&[-2,5]\\
s899&$b$&4&s900&$b=Cx;c$&[4,-1]\\
s938&$c$&0=[-3]&s939&$c$&0=[-3]\\
\hline
\end{tabular}
\end{table}

\begin{table}
\centering
\begin{tabular}{|lrr|lrr|}
\hline
Name&Standard&Poly&Name&Standard&Poly\\
\hline
v0895&$b$&[2,-14]&v0896&$a$&2\\
v0948&$b=Ax,c=Ay;a$&[-5,2]&v0949&$b=Cx,a=yC^2;c$&[5,2]\\
v0950&$b$&-3&v0951&$c$&3\\
v1000&$b$&3&v1001&$b=Ax;a$&-3\\
v1016&$b=xa;a$&4&v1017&$c$&-4\\
v1066&$b=xc;c$&5&v1067&$a$&-5\\
v1083&$a$&-5&v1084&$a=xb;b$&-5\\
v1095&$b$&[7,-4]&v1096&$b=ax;a$&[-7,-4]\\
v1097&$b=xA;a$&[7,6]&v1098&$c=Ax;a$&[7,-6]\\
v1104&$b=Cy,c=az;a$&[-7,8]&v1105&$c=xA;a$&[-7,-8]\\
v1110&$c=xA;a$&7&v1111&$b=xA;a$&-7\\
v1123&$c=xa;a$&[-8,-6]&v1124&$c=Ax;a$&[8,-6]\\
v1128&$c=A^2x;a$&[8,-10]&v1129&$c=Ax;a$&[-8,-10]\\
v1491&$c$&-2&v1492&$b$&-2\\
v1684&$c=Ax,a=yb^2;b$&[2,-2,-3]&v1737&$a=yz^2yz^3,b=Z^2Y;z$&[-4,4,-3]\\
v1781&$a$&-3&v1782&$c$&3\\
v1793&$b=Ax,a=yz^2,$&&v1858&$a=xC^2;c$&(3,-1)\\
&$c=Z^3Y;z$&[4,-4,5]&&&\\
v1863&$a=C^2x;c$&[2,-2,7]&v1893&$a=B^2x;b$&[-3,3,2]\\
v1897&$c$&-2&v1898&$c$&2\\
v1901&$a=xB;b$&[-4,1]&v1902&$a=xB;b$&[4,1]\\
v2001&$b=Cx;c$&[4,-7]&v2002&$b=xC;c$&[4,7]\\
v2022&$b$&[-3,-15]&v2023&$c$&-3\\
v2037&$c=A^3x;a$&(3,1)&v2066&$a=xy^3,b=Y^2X;y$&[3,-3,8]\\
v2103&$b=C^2x,c=X^2y;x$&[-5,5,-3]&v2130&$c=az,a=xy^2,$&\\
&&&&$b=Y^3X;y$&[-5,5,-2]\\
v2134&$b$&3&v2135&$c$&-3\\
v2146&$c$&4&v2147&$a$&4\\
v2151&$a=xy^3,c=Y^2X;y$&[-5,5,-7]&v2174&$b=Ax;a$&[-5,-1]\\
v2175&$c$&-5&v2182&$b$&5\\
v2183&$b$&-5&v2205&$a=xy^2,c=Y^3X;y$&[-5,5,-8]\\
v2257&$c=Bx;b$&[-5,11]&v2258&$a$&-5\\
v2304&$c=Bx;b$&[5,-14]&v2305&$a$&5\\
v2308&$c=xa,b=yX^2$,&&v2346&$c$&3\\
&$a=zx^3;x$&[-3,3,-3,2]&&&\\
v2347&$a$&3&v2365&$a=xy^3,c=Y^2X;y$&[-3,3,-3,4]\\
v2388&$b$&-2&v2389&$a$&2\\
v2438&$b=A^3x;a$&[2,-2,0,3]&v2467&$c$&-2\\
v2468&$b$&2&{\bf v2530}&$a=bx,c=yX^2,$&\\
&&&&$b=xz;x$&({\bf 1})\\
\hline
\end{tabular}
\end{table}

\begin{table}
\centering
\begin{tabular}{|lrr|lrr|}
\hline
Name&Standard&Poly&Name&Standard&Poly\\
\hline
v2575&$c$&3&v2576&$a$&[-3,10]\\
v2605&$c=xA^2;a$&[-3,-1,0]&v2706&$c=Ax;a$&[-3,-14]\\
v2707&$a$&3&v2708&$c=Ax;a$&[-3,1,0]\\
v2743&$a$&4&v2744&$a$&[-4,-7]\\
v2787&$a$&(-2,0)&v2807&$a$&6\\
v2808&$a$&[-6,6]&v2861&$a$&-2\\
v2862&$b$&2&{\bf v2869}&$c=Ax,b=Ay,$&\\
&&&&$a=zx^2;x$&[{\bf -1},2,-2,3]\\
v2874&$b$&(-2,3)&v2926&$c=Ax;a$&[6,0]\\
v2927&$c=Ax;a$&[-6,0]&v2997&$a$&2\\
v2998&$b$&-2&v3003&$b=Ax;a$&[-6,-3]\\
v3004&$a$&-6&v3036&*2 gen* $a$&[3,4,5]\\
v3092&$c=xa^3,b=Ay;a$&[2,-2,0,2,-1]&{\bf v3093}&$a=xB^2,c=yB;
b$&[{\bf -1},1,1]\\
v3102&$b$&4&v3103&$b$&-4\\
v3145&$c=By,a=bx;b$&[2,1,2]&v3168&$a$&-3\\
v3169&$b$&-3&v3188&$c$&-2\\
v3189&$a$&2&v3210&$c$&2\\
v3219&$a$&0=[4]&v3221&$c$&0=[4]\\
v3226&$b$&-2&v3228&$c$&-4\\
v3243&$a=B^2x,c=b^2y;b$&[-2,-1,-2]&v3244&$b=A^2x;a$&[-2,1,-2]\\
v3245&$a=C^3x,b=c^4y;c$&[2,-1,-1,3]&v3272&$c$&[3,-10]\\
v3273&$a$&3&v3293&$b=xA;a$&(-2)\\
v3329&$a=b^3x,c=By;b$&[2,-1,-1,2,-1]&v3337&$b$&[-4,-10]\\
v3338&$b$&4&v3377&$b=xA;a$&(-2,-1)\\
v3382&$c$&-5&v3394&$b$&[-3]\\
v3395&$b=Ax;a$&[3]&v3452&$c$&[4,-2]\\
v3453&$a$&[-4,-2]&v3492&$c$&-2\\
v3493&$a$&2&v3498&$a=C^3x,b=cy;c$&[2,0,-2,3]\\
v3526&$b=Ax,c=a^2y;a$&(-2)&{\bf v3541}&$b=ax,a=yz^4,$&\\
&&&&$c=Z^3Y;z$&[{\bf 1},-2,1,0,2,-3]\\
\hline
\end{tabular}
\end{table}
\begin{table}
\caption{Cusped virtually fibred census 3-manifolds}
\centering
\begin{tabular}{|lcl|lcl|}
\hline
Name&Ratio&Name of fibred&Name&Ratio&Name of fibred\\
\hline
m006,m007&1/2&v1241&
m015,m017&1/1&m015=5a1($5_2=3 2$)\\
m029,m030&1/2&v$3140^{\bf 2}$&
m032,m033&1/1&m032=6a3($6_1=4 2$)\\
m035,m037&1/1&m039,m040&
m045,m046&1/2&v$3383^{\bf 3}$\\
\y&1/2&v$3218^{\bf 2}$,v$3220^{\bf 2}$,v$3222^{\bf 2}$,&&&\\
&&v$3225^{\bf 2}$,v$3227^{\bf 3}$&&&\\
m053,m054&1/1&m053=7a4($7_2=5 2$)&
m073,m074&1/1&m074=8a11($8_1=6 2$)\\
m079,m080&1/2&10a$171^{\bf 2}(2 6 2)$&
m093,m094&1/1&m094=9a27($9_2=7 2$)\\
m139&&$\mathcal A 1$&
m148,m149&1/2&8a$31^{\bf 2}$($8^2_4=3 2 3$)\\
m208&&$\mathcal A 3$&
m287,m288&1/2&9a39($9_{10}=3 3 3$)\\
m306,m307&1/1&s298,s299&
m340&1/1&m340=7a5($7_3=4 3$)\\
m410&&$\mathcal A 3$&&&\\
\hline
s016,s017&1/1&s016=10a75($10_1=8 2$)&
s022,s023&1/1&s023=11a247($9 2$)\\
s119&&$\mathcal A 3$&
s348&1/1&m$329^{\bf 2}$\\
s349,s350&1/1&m$328^{\bf 2}$&
s423,s424&1/1&m$359^{\bf 2}$\\
s437&1/1&m$367^{\bf 2}$&
s477&1/1&m$391^{\bf 2}$\\
s478,s480&1/1&s479,v0953&
s558&1/1&s588=9a38($9_3=6 3$)\\
s643,s644&1/2&11a333($4 1 1 1 4$)&
s648,s649&1/1&v1241,\\
&&&&&s648=7a6($7_4$=313)\\
s673,s674&1/1&v1276,v1277&
s725,s726&1/1&s726=8a18($8_3=4 4$)\\
s763,s764&1/2&9a16($9_{23}=2 2 1 2 2$)&
s772,s773,&&$\mathcal A 7$\\
&&&s779,s784&&\\
s788&1/1&s789,v1539,v1540&
s818,s819&1/1&s817,v1638\\
s862&1/1&s862=8a17($8_4=4 1 3$)&
s870&1/1&s870=9a35($9_4=5 4$)\\
s899,s900&2/1&m015=5a1($5_2=3 2$)&&&\\
\hline
\end{tabular}
\end{table}

\begin{table}
\centering
\begin{tabular}{|lcl|lcl|}
\hline
Name&Ratio&Name of fibred&Name&Ratio&Name of fibred\\
\hline
v0016,v0017&1/1&v0016=12a803&
v0024,v0025&1/1&v0025=13a3143\\
&&($[10]2$)&&&($[11]2$)\\
v0571&1/1&m340=7a5($7_3=4 3$)&
v0785&1/1&m$357^{\bf 2}$\\
v0819&1/1&m$366^{\bf 2}$,v0820&
v0954&1/1&m$388^{\bf 2}$\\
v1010,v1012&1/1&s$506^{\bf 2}$&
v1011&1/1&s$503^{\bf 2}$\\
v1035,v1036&1/2&11a365($353$)&
v1112,v1113&1/1&s$549^{\bf 2}$\\
v1152&1/1&s$568^{\bf 2}$&
v1168&1/1&s$577^{\bf 2}$\\
v1172&1/1&s$578^{\bf 2}$&
v1179&1/1&v$1178^{\bf 2}$\\
v1194&1/1&s$602^{\bf 2}$&
v1205&1/1&v$1204^{\bf 2}$\\
v1210&1/1&s$621^{\bf 2}$&
v1229&1/1&s$638^{\bf 2}$\\
v1243&1/1&v1243=11a364($83$)&
v1256&1/1&s$661^{\bf 2}$\\
v1676&1/1&s$831^{\bf 2}$&
v1858&&$\mathcal A 1$\\
v2018&1/1&s$876^{\bf 2}$&
v2037&1/1&s$880^{\bf 2}$\\
v2078&1/1&s$887^{\bf 2}$&
v2158&1/1&s$895^{\bf 2}$\\
v2203&1/1&s$898^{\bf 2}$,v2202&
v2238&1/1&s$906^{\bf 2}$\\
v2284,v2285&1/1&v2284=9a36($9_5=513$)&
v2297,v2298&1/1&v2296\\
v2339&1/1&s$914^{\bf 2}$&
v2346,v2347&1/1&v2345\\
v2361,v2362&1/1&v2362=10a117&
v2467,v2468&1/1&v2469\\
&&($10_3=64$)&&&\\
v2488&1/1&v2488=10a113&
v2520&1/1&v2520=11a342\\
&&($10_4=613$)&&&($74$)\\
v2575,v2576&1/1&v2574&
v2706,v2707&1/1&v2705\\
v2787&&$\mathcal A 2$&
v2796,v2797&1/2&11a119($23132$)\\
v2858&1/1&v2858=10a114&v2874&&$\mathcal A 3$\\
&&($10_8=514$)&&&\\
v2894&1/1&11a358($65$)&
v$2943^{\bf 2}$,v2944&1/1&v$2942^{\bf 2}$\\
v3128&1/1&v3126,v$3127^{\bf 2}$&
v3210&1/1&v3207,v3208,\\
&&&&&v3209\\
v3243,v3244&1/1&v3246,v3247&
v3310&1/1&v3310=7a3\\
&&&&&($7_5=322$)\\
v3377&1/1&v$3376^{\bf 2}$,v3378&
v$3379^{\bf 2}$,v$3384^{\bf 2}$&1/1&v$3383^{\bf 3}$\\
v$3396^{\bf 2}$&1/1&v$3393^{\bf 2}$&
v3427&1/1&v$3426^{\bf 2}$\\
v3457&1/1&v$3456^{\bf 2}$&
v3492,v3493&1/1&v3490,v3491\\
\hline
\end{tabular} 
\end{table}

\begin{table}
\caption{Closed fibred census 3-manifolds}
\centering
\begin{tabular}{|llrcl|} \hline 
Volume & Name & $\mathbb Z$ & neg & Alternative\\

\hline
3.1663 &m160(3,1)&&&\\
\y&m159(4,1)&&& \\
3.1772 &m199(-4,1)&2&&\\
\y&m122(-4,1)&&&\\
3.6638 &s942(-2,1)&&-&s957(-1,2)\\
\y&m336(-1,3)&&&\\
3.7028&m345(1,2)&\s&&\\
3.7708&m289(7,1)&2&&\\
\y&m280(1,4)&&&\\
3.8534&m304(-5,1)&&&\\
\y&m305(-1,3)&&&\\
3.9466&s385(5,1)&3&&\\
3.9702&s296(-1,3)&&&\\
\y&s297(5,1)&&&\\
4.0597&s912(0,1)&2&&\\
\y&m401(-2,3)&&&\\
\y&m371(-1,3)&&&\\
\y&m368(-4,1)&&&\\
4.4081&s580(-5,1)&2&&\\
\y&s581(-1,3)&&&\\
4.4153&s869(-1,2)&\s&&\\
\y&s861(3,1)&&&\\
4.4191&v1191(-5,1)&&&\\
\y&v1076(-5,1)&&&\\
4.4646&s924(3,1)&&&\\
\y&v1408(4,1)&&&\\
4.5169&s677(1,3)&&&\\
\y&s676(5,1)&&&\\
\hline
\end{tabular}
\end{table}
\begin{table}
\centering
\begin{tabular}{|llrcl|} \hline 

Volume & Name & $\mathbb Z$ & neg & Alternative\\

\hline
4.5559&v2641(-4,1)&\s&&\\
\y&s745(3,2)&&&\\
4.6307&s646(5,2)&&&\\
4.7135&v1539(5,1)&$\beta$2&&\\
\y&s789(-5,1)&&&v1540(1,3)\\
4.7252&s719(7,1)&&&\\
\y&v1373(-2,3)&&&\\
4.7517&v3209(3,1)&&&v3514(-2,1)\\
\y&v2420(-3,1)&&&\\
4.7659&v2099(-4,1)&&&\\
\y&v2101(3,1)&&&\\
4.7740&s789(5,1)&&&v1670(-1,3)\\
\y&v1539(-5,1)&&&\\
4.7874&v1721(1,4)&\s&&\\
4.9068&v2771(-4,1)&&&\\
4.9069&s836(-6,1)&4&&\\
4.9094&v2986(1,2)&\s&&\\
4.9717&v2209(2,3)&\s&&\\
5.1171&v2054(-7,1)&&&\\
5.1379&v3066(-1,2)&\s&&\\
\y&v2563(5,1)&&&\\
\y&v2345(5,1)&&&\\
5.1706&v3209(-3,1)&&&v3486(3,1)\\
5.1984&v3077(5,1)&&-&\\
\y&v2959(-3,1)&&-&\\
5.2007&v2671(-2,3)&\s&&\\
5.2983&s928(2,3)&\s&&\\
\hline
\end{tabular}
\end{table}
\begin{table}
\centering
\begin{tabular}{|llrcl|} \hline 

Volume & Name & $\mathbb Z$ & neg & Alternative\\

\hline
5.3334&v3390(3,1)&&&\\
\y&v3209(4,1)&&&\\
\y&v2913(-3,2)&&&\\
\y&v3505(-3,1)&2&&\\
\y&v3261(4,1)&&&\\
\y&v3262(3,1)&&&\\
5.3488&v2678(-5,1)&&&\\
5.4633&v3107(3,2)&\s&&\\
5.4957&v3216(4,1)&&&\\
\y&v3217(-1,3)&&&\\
5.4962&v3320(4,1)&3&&\\
5.5410&v3091(-2,3)&&&\\
5.5636&v3214(1,3)&&&\\
\y&v3215(-4,1)&&&\\
5.5736&v3209(-4,1)&&&\\
5.6510&v2984(-1,3)&\s&&\\
5.6664&v3209(5,1)&&&\\
5.6743&v3019(5,2)&\s&&\\
5.7024&v3212(1,3)&\s&&\\
5.8111&v3209(-5,1)&&&\\
5.8524&v3425(-3,2)&&&\\
5.8664&v3209(6,1)&&&\\
5.8760&v3318(4,1)&\s&&\\
5.9780&v3352(1,4)&\s&&\\
6.0075&v3398(2,3)&\s&&\\
6.0502&v3378(-1,4)&\s&&\\
6.1102&v3408(1,3)&\s&&\\
6.1203&v3467(-2,3)&\s&&\\
6.1254&v3445(6,1)&\s&&\\
6.2391&v3509(4,3)&&&\\
\y&v3508(4,1)&&&\\
6.2428&v3504(-2,3)&&&\\
\hline
\end{tabular}
\end{table}
\begin{table}
\caption{Closed non-fibred census 3-manifolds}
\centering
\begin{tabular}{|llrr|}
\hline
Volume&Name&Standard&Poly\\
\hline
4.4559&s528(-1,3)&$a=xB;b$&[-2,2,1]\\
\y&s527(-5,1)&$a=xB;b$&[2,2,-1]\\
4.5760&s644(-4,3)&$a=xB^2;b$&[2,-2,5]\\
\y&s643(-5,1)&$a=xB;b$&[2,2,5]\\
4.7494&{\bf v2018(-4,1)}&$a=xB;b$&(2,-1,1)\\
4.7809&v1436(-5,1)&$b$&[3,-1,3]\\
4.7904&s750(4,3)&$a=xy^3,b=Y^2X;y$&[-3,3,5]\\
\y&s749(5,1)&$a=xB;b$&[3,3,5]\\
4.8461&s789(-5,2)&&\texttt x\\
\y&v1539(5,2)&&\texttt x\\
4.8511&{\bf v2238(-5,1)}&$a=xB;b$&(2,1,-1)\\
\y&v3209(1,2)&&\texttt x\\
\y&s828(-4,3)&$a=xB;b$&[2,4,5]\\
4.8810&v1695(5,1)&$a$&[3,-2,3]\\
5.0362&s862(7,1)&(The knot $8_4)\,b$&[-2,5,-5]\\
\y&v2190(4,1)&$a=xB;b$&[2,5,5]\\
5.2283&v3209(-1,2)&&\texttt x\\
\y&v2593(4,1)&$a=xB;b$&[2,4,3]\\
5.3811&v3209(3,2)&&\texttt x\\
\y&v3027(-3,1)&$a=xB;b$&[2,4,7]\\
5.4334&v2896(-6,1)&$a=xB;b$&[-2,3,0]\\
\y&v2683(-6,1)&$a=xB;b$&[2,3,0]\\
5.4561&v2796(4,1)&$a$&[2,-1,5]\\
\y&v2797(-3,4)&$a$&[2,1,5]\\
5.5573&v2948(-6,1)&$a=xB;b$&[3,-2,0]\\
\y&v2794(-6,1)&$a=xB;b$&[3,2,0]\\
5.5736&{\bf v3183(-3,2)}&$a=xy^3,b=Y^2X;y$&[-2,0,0]\\
5.6562&v3145(3,2)&$b$&[-2,-1,-2]\\
\y&v3181(-3,2)&$a=xB;b$&[2,5,8]\\
5.6872&v3036(3,2)&$b$&[3,4,5]\\
5.7024&v3209(1,3)&&\texttt x\\
\y&v3269(4,1)&$a$&[3,6,7]\\
5.7057&v3209(-3,2)&&\texttt x\\
5.7243&v3209(2,3)&&\texttt x\\
\y&v3313(3,1)&$b$&[3,6,8]\\
5.8041&v3239(3,2)&$a=xB;b$&[3,5,7]\\
5.8060&v3209(5,2)&&\texttt x\\
5.8073&v3209(-1,3)&&\texttt x\\
5.8759&v3209(4,3)&&\texttt x\\
5.8882&v3244(4,3)&$C=ax,a=yz^2,b=Z^3Y;z$&[2,-1,2]\\
\y&v3243(-4,1)&$a=xc,b=yx,c=Z^2Y;z$&[2,1,2]\\
\hline
\end{tabular}
\end{table}
\begin{table}
\caption{Closed virtually fibred census 3-manifolds}
\centering
\begin{tabular}{|llcll|}
\hline
Volume&Name&neg&Ratio&Name of fibred\\
\hline
0.9427&m003(-3,1)&&1/4
&m289(7,1)\\
&&&&m280(1,4)\\&&&&\\
1.0149&m007(3,1)&-&1/1&m$007(3,1)^0$\\&&&&\\
1.4140&m009(4,1)&&3/8
&m289(7,1)\\
&&&&m280(1,4)\\&&&&\\
1.5831&m007(4,1)&&1/2
&m160(3,1)\\
&&&&m159(4,1)\\&&&&\\
1.5886&m006(3,1)&-&1/3&v2099(-4,1)\\
&m003(-5,4)
&&&v2101(3,1)\\&&&&\\
1.8319&m009(5,1)&&1/3&v3217(-1,3)\\
&m010(-2,3)&&&\\
&m009(-5,1)&&&\\
&m006(1,3)&&&\\
&&&&\\
1.8854&m007(5,1)&&1/2&m289(7,1)\\
&m006(-1,3)&&&m280(1,4)\\&&&&\\
2.0298&m036(-3,2)&-&2/1&m$007(3,1)^0$\\
&m010(-4,3)&&&\\&&&&\\
\y&m010(4,1)&&1/2&m371(-1,3)\\
&&&&m368(-4,1)\\&&&&\\
2.5689&m039(6,1)&&2/3&m304(-5,1)\\
&m035(-6,1)&&&m305(-1,3)\\
&m037(2,3)&&&\\
&m130(-3,1)?&&&\\
&m120(-4,1)&-&&\\
&m223(3,1)&&&\\
&m038(-6,1)&&&\\
&m036(-2,3)&-&&\\
&&&&\\
\hline
\end{tabular} 
\end{table}
\begin{table}
\centering
\begin{tabular}{|llcll|}
\hline
Volume&Name&neg&Ratio&Name of fibred\\
\hline
2.6667&m135(-1,3)&&1/2&v3505(-3,1)\\
&m135(1,3)&-&&v3261(4,1)\\
&m168(3,2)&&&v3262(3,1)\\
&m140(4,1)&-&&\\
&&&&\\
2.8281&m221(3,1)&&3/4&m289(7,1)\\
&m070(1,4)&&&m280(1,4)\\
&m139(2,3)&&&\\
&&&&\\
3.0448&m247(-1,3)&-&3/1&m$007(3,1)^0$\\&&&&\\
3.1772&m303(-3,1)&&1/1&m199(-4,1)\\
&&&&m122(-4,1)\\&&&&\\
\y&m141(4,1)&&2/3&v2099(-4,1)\\
&m249(1,2)&&&v2101(3,1)\\
&s254(-3,1)&&&\\
&s479(-3,1)&-&&\\
&m146(-2,3)&&&\\
&m188(4,1)&&&\\
&m148(6,1)&&&\\
&m149(-2,3)&&&\\
&m206(3,2)&&&\\
&m159(-2,3)&&&\\
&&&&\\
3.6638&s960(-1,2)&&1/1&s942(-2,1)\\
&m304(5,1)&&&m336(-1,3)\\
&&&&\\
\y&s572(1,2)&-&2/3&v3216(4,1)\\
&m293(4,1)&&&\\
&s645(-1,2)&&&\\
&s297(-1,3)&&&\\
&s778(-3,1)&&&\\
&s775(-1,2)&&&\\
&s682(-3,1)&&&\\
&&&&\\
\hline
\end{tabular} 
\end{table}
\begin{table}
\centering
\begin{tabular}{|llcll|}
\hline
Volume&Name&neg&Ratio&Name of fibred\\
\hline
\y&s296(5,1)&&2/3&v3217(-1,3)\\
&s779(2,1)&-&&\\
&m312(-1,3)&&&\\
&s595(3,1)&&&\\
&s775(-3,1)&&&\\
&s350(-4,1)&&&\\
&m294(4,1)&&&\\
&s495(1,2)&-&&\\
&&&&\\
3.7708&m369(-3,2)&&1/1&m289(7,1)\\
&m371(3,2)&&&m280(1,4)\\
&s478(-1,2)&&&\\
&s479(1,2)&-&&\\
&&&&\\
3.9702&s784(1,2)&-&1/1&s296(-1,3)\\
&m303(5,1)&&&s297(5,1)\\
&m376(3,2)&&&\\
&&&&\\
4.0597&v0825(4,1)&&4/1&m$007(3,1)^0$\\
&m358(1,3)&&&\\
&s775(1,2)&&&\\
&s778(-3,2)&&&\\
&s779(1,2)&-&&\\
&m395(-2,3)&&&\\
&s787(1,2)&&&\\
&s440(-1,3)&&&\\&&&&\\
\y&s705(-3,1)&&1/1&m371(-1,3)\\
&s779(-3,2)&&&m368(-4,1)\\
&s772(-3,2)&&&\\
&&&&\\
4.2421&v2101(1,2)&&5/4&m289(7,1)\\
&&&&m280(1,4)\\&&&&\\
\hline
\end{tabular} 
\end{table}
\begin{table}
\centering
\begin{tabular}{|llcll|}
\hline
Volume&Name&neg&Ratio&Name of fibred\\
\hline

4.4153&v2101(-1,3)&&1/1&s869(-1,2)\\
&s779(-4,1)&&&s861(3,1)\\
&s775(-4,1)&&&\\
&s778(3,1)&&&\\
&s772(-4,1)&&&\\
&s773(3,1)&&&\\
&s786(3,1)&&&\\
&s781(-4,1)&&&\\
&&&&\\
4.4646&s781(-2,3)&&1/1&s924(3,1)\\
&s786(-1,3)&&&v1408(4,1)\\
&s773(-1,3)&-&&\\
&s777(-5,1)&&&\\
&v2787(1,2)&-&&\\
&&&&\\
4.6307&s645(5,2)&&1/1
&s646(5,2)\\&&&&\\
4.7135&s889(3,2)&&5/4&m289(7,1)\\
&v2739(1,2)&&&m280(1,4)\\
&&&&\\
\y&v2797(2,1)&&1/1&v1539(5,1)\\
&v2573(-3,2)&&&s789(-5,1)\\
&s788(-1,3)&&&\\
&&&&\\
4.7494&{\bf v2018(-4,1)} $\beta 1$&&3/2&m160(3,1)\\
&&&&m159(4,1)\\
4.7659&v2787(-3,1)&&3/2&m199(-4,1)\\
&v1644(-2,3)&&&m122(-4,1)\\
&v2100(-3,1)&&&\\
&&&&\\
\y&s916(-3,2)&&1/1&v2099(-4,1)\\
&s957(1,2)&&&v2101(3,1)\\
&s821(2,3)&&&\\
&s960(1,2)&&&\\
&&&&\\
\hline
\end{tabular}
\end{table}
\begin{table}
\centering
\begin{tabular}{|llcll|}
\hline
Volume&Name&neg&Ratio&Name of fibred\\
\hline
4.9068&v2018(2,3)&-&1/1
&v2771(-4,1)\\
&&&&\\
5.0747&v3216(-4,1)&&5/1&m$007(3,1)^0$\\
&v3210(3,1)&&&\\
&v2636(2,3)&&&\\
&v2417(-1,3)&&&\\&&&&\\
\y&v3213(-3,1)&&5/4
&m371(-1,3)\\
&&&&m368(-4,1)\\&&&&\\
5.1379&v3100(-3,1)&&4/3&m304(-5,1)\\
&v2346(-1,3)&&&m305(-1,3)\\
&v2345(-1,3)&&&\\
&v3469(3,1)&&&\\
&v3106(1,3)&&&\\
&s916(5,1)&&&\\
&v3214(-3,1)&&&\\
&&&&\\
\y&v2346(5,1)&&1/1
&v2563(5,1)\\
&&&&v2345(5,1)\\&&&&\\
5.3334&v3210(-3,1)&&&v3209(4,1)\\
&v3207(-3,1)&&&v3505(-3,1)\\
&v3208(4,1)&&&v3261(4,1)\\
&v3106(-3,1)&&&v3262(3,1)\\
&v3107(-4,1)&&&\\
&v3331(-2,3)&&&\\
&&&&\\
5.4957&v3213(-1,3)&&1/1
&v3216(4,1)\\&&&&\\
\y&v3412(5,1)&&1/1
&v3217(-1,3)\\&&&&\\
5.6562&v3387(3,2)&&3/2&m289(7,1)\\
&v3136(-1,3)&&&m280(1,4)\\
&&&&\\
\hline
\end{tabular}
\end{table}

\end{document}